\documentclass[10pt]{article}
\usepackage{amsfonts}

\begin {document}

\topmargin= -.2in \baselineskip=20pt
\title {$L$-Functions for Symmetric Products of Kloosterman Sums}
\author {Lei Fu\\
{\small Institute of Mathematics, Nankai University, Tianjin, P.
R. China}\\
{\small leifu@nankai.edu.cn}\\{}\\
Daqing Wan\\
{\small Institute of Mathematics, Chinese Academy of Sciences,
Beijing, P. R. China}
\\
{\small Department of Mathematics, University of California,
Irvine, CA
92697}\\
{\small dwan@math.uci.edu}}
\date{}
\maketitle
\bigskip
\centerline {\bf Abstract}
\bigskip

The classical Kloosterman sums give rise to a Galois
representation of the function field unramified outside $0$ and
$\infty$. We study the local monodromy of this representation at
$\infty$ using $l$-adic method based on the work of Deligne and
Katz. As an application, we determine the degrees and the bad
factors of the $L$-functions of the symmetric products of the
above representation. Our results generalize some results of Robba
obtained through $p$-adic method.


\bigskip
\centerline {\bf 0. Introduction}

\bigskip
\bigskip
Let ${\bf F}_q$ be a finite field of characteristic $p$ with $q$
elements and let $\psi: {\bf F}_q\to \overline {\bf Q}_l^\ast$ be
a nontrivial additive character. Fix an algebraic closure ${\bf
F}$ of ${\bf F}_q$. For any integer $k$, let ${\bf F}_{q^k}$ be
the extension of ${\bf F}_q$ in ${\bf F}$ with degree $k$. If
$\lambda$ lies in ${\bf F}_{q^k}$, we define the $(n-1)$-variable
Kloosterman sum by
$${\rm Kl}_n({\bf F}_{q^k}, \lambda)=\sum\limits_{x_1\cdots x_n=\lambda
,\;x_i\in  {\bf F}_{q^k}} \psi({\rm Tr}_{{\bf F}_{q^k}/{\bf
F}_q}(x_1+\cdots + x_n )).$$ For any $\lambda\in {\bf F}^\ast $,
define ${\rm deg}(\lambda)=[{\bf F}_q(\lambda):{\bf F}_q]$. The
$L$-function $L(\lambda, T)$ associated to Kloosterman sums is
defined by
$$L(\lambda, T)=\exp \biggl(\sum\limits_{m=1}^\infty
{\rm Kl}_n({\bf F}_{q^{m{\rm deg}(\lambda)}},
\lambda)\frac{T^m}{m}\biggr).$$ One can show that
$L(\lambda,T)^{(-1)^n}$ is a polynomial of degree $n$ with
coefficients in ${\bf Z}[\zeta_p]$, where $\zeta_p$ is a primitive
$p$-th root of unity. This follows from 7.4 and 7.5 in [D1]. Write
\begin{eqnarray*}
L(\lambda, T)^{(-1)^n}&=&(1-\pi_1(\lambda)T)\cdots
(1-\pi_n(\lambda)T).
\end{eqnarray*}
Then for any positive integer $m$, we have
$${\rm Kl}_n({\bf F}_{q^{m{\rm deg}(\lambda)}}, \lambda)
=(-1)^{n-1}(\pi_1(\lambda)^m+\cdots+ \pi_n(\lambda)^m).$$

We have a family of $L$-functions $L(\lambda, T)$ parameterized by
the parameter $\lambda$. Let $|{\bf G}_m|$ be the set of Zariski
closed points in ${\bf G}_m={\bf P}^1-\{0,\infty\}$. This is the
parameter space for $\lambda$. For a positive integer $k$, the
$L$-function for the $k$-th symmetric product of the Kloosterman
sums is defined by
$$L({\bf G}_m, {\rm Sym}^k({\rm Kl}_n), T) = \prod_{\lambda\in
|{\bf G}_m|} \prod_{i_1+\cdots +i_n=k}
(1-\pi_1(\lambda)^{i_1}\cdots \pi_n(\lambda)^{i_n} T^{{\rm
deg}(\lambda)})^{-1}.
$$
This is a rational function in $T$ by Grothendieck's formula for
$L$-functions, see 3.1 of [Rapport] in [SGA 4$\frac{1}{2}$]. In
the proof of Lemma 2.2, we shall see that $L({\bf G}_m, {\rm
Sym}^k({\rm Kl}_n), T)$ actually has coefficients in ${\bf Z}$.
Thus, this rational function is geometric in nature, that is, it
should come from the zeta function of some varieties or motives.
What are these varieties and motives? From arithmetic point of
view, our fundamental problem here is to understand this sequence
of $L$-functions with integer coefficients parameterized by the
arithmetic parameter $k$ and its variation with $k$, from both
complex as well as $p$-adic point of view. There is still a long
way to go toward a satisfactory answer of this basic question,
most notably from $p$-adic point of view. We shall make some
remarks at the end of this introduction section.

In [R], Robba studied the $L$-function $L({\bf G}_m, {\rm
Sym}^k({\rm Kl}_2), T)$ in the special case $n=2$ using Dwork's
$p$-adic methods. He conjectured a degree formula, obtained the
functional equation and the bad factors of the $L$-function
$L({\bf G}_m, {\rm Sym}^k({\rm Kl}_2), T)$. On the other hand,
Kloosterman sums have also been studied in great depth by Deligne
and Katz using $\l$-adic methods. The purpose of this paper is to
use their fundamental results to derive as much arithmetic
information as we can about this sequence of $L$-functions $L({\bf
G}_m, {\rm Sym}^k({\rm Kl}_n), T)$ for $n>1$ not divisible by $p$.
We have

\bigskip
\noindent {\bf Theorem 0.1.} Suppose $(n,p)=1$. Let $\zeta$ be a
primitive $n$-th root of unity in ${\bf F}$. For a positive
integer $k$,  let $d_k(n,p)$ denote the number of $n$-tuples
$(j_0,j_1,\ldots, j_{n-1})$ of non-negative integers satisfying
$j_0+j_1+\cdots + j_{n-1}=k$ and $j_0+j_1\zeta+\cdots + j_{n-1}
\zeta^{n-1}=0$ in ${\bf F}$. Then, the degree of the rational
function $L({\bf G}_m, {\rm Sym}^k({\rm Kl}_n), T)$ is
$$\frac{1}{n}\left({k+n-1\choose n-1} - d_k(n,p)\right).$$
\bigskip

Note that in most cases, the $L$-function $L({\bf G}_m, {\rm
Sym}^k({\rm Kl}_n), T)$ is a polynomial, not just a rational
function. This is the case for instance if $pn$ is odd or $n=2$.
This follows from Katz's global monodromy theorem for the
Kloosterman sheaf. Recall that the Kloosterman sheaf is a lisse
$\overline{\bf Q}_l$-sheaf ${\rm Kl}_n$ on ${\bf G}_m$ such that
for any $x\in {\bf G}_m({\bf F}_{q^k})={\bf F}_{q^k}^\ast$, we
have
$${\rm Tr}(F_x, {\rm Kl}_{n,\bar x})=(-1)^{n-1} {\rm Kl}_n({\bf F}_{q^k},
x),$$ where $F_x$ is the geometric Frobenius element at the point
$x$. For any Zariski closed point $\lambda$ in ${\bf G}_m$, let
$\pi_1(\lambda),\ldots, \pi_n(\lambda)$ be all the eigenvalues of
the geometric Frobenius element $F_\lambda$ on ${\rm Kl}_{n,\bar
\lambda}$. Then
\begin{eqnarray*}
{\rm Kl}_n({\bf F}_{q^{m{\rm
deg}(\lambda)}},\lambda)&=&(-1)^{n-1}{\rm Tr}(F_\lambda^m, {\rm
Kl}_{n,\bar \lambda})\\
&=& (-1)^{n-1}(\pi_1(\lambda)^m+\cdots +\pi_n(\lambda)^m).
\end{eqnarray*}
So we have
\begin{eqnarray*}
L(\lambda, T)^{(-1)^n}&=&\exp \biggl((-1)^n\sum_{m=1}^\infty {\rm
Kl}_n({\bf
F}_{q^{m{\rm deg}(\lambda)}},\lambda)\frac{T^m}{m}\biggr)\\
&=& \exp\biggl(\sum_{m=1}^\infty
(-(\pi_1(\lambda)^m+\cdots+\pi_n(\lambda)^m))\frac{T^m}{m}\biggr)\\
&=& (1-\pi_1(\lambda)T)\cdots(1-\pi_n(\lambda)T).
\end{eqnarray*}
We have
$$\prod_{i_1+\cdots+i_n=k}(1-\pi_1(\lambda)^{i_1}\cdots
\pi_n(\lambda)^{i_n}T^{{\rm deg}(\lambda)})= {\rm det}(1-F_\lambda
T^{{\rm deg}(\lambda)}, {\rm Sym}^k ({\rm Kl}_{n,\bar\lambda})).$$
Therefore the $L$-function for the $k$-th symmetric product of the
Kloosterman sums $$L({\bf G}_m,{\rm Sym}^k({\rm Kl}_n),
T)=\prod_{\lambda \in |{\bf
G}_m|}\prod_{i_1+\cdots+i_n=k}(1-\pi_1(\lambda)^{i_1}\cdots
\pi_n(\lambda)^{i_n}T^{{\rm deg}(\lambda)})^{-1}$$ is nothing but
the Grothendieck $L$-function of the $k$-th symmetric product of
the Kloosterman sheaf:
$$L({\bf G}_m,{\rm Sym}^k({\rm Kl}_n), T)=\prod_{\lambda \in |{\bf
G}_m|} {\rm det}(1-F_\lambda T^{{\rm deg}(\lambda)}, {\rm Sym}^k
({\rm Kl}_{n,\bar\lambda}))^{-1}.$$

The Kloosterman sheaf ramifies at the two points $\{0, \infty\}$.
We would like to determine the bad factors of the $L$-function
$L({\bf G}_m, {\rm Sym}^k({\rm Kl}_n), T)$ at the two ramified
points. Assume $n|(q-1)$. Then we can explicitly determine the bad
factor at $\infty$. The key is to determine the local monodromy of
the Kloosterman sheaf at $\infty$, that is, to determine ${\rm
Kl}_n$ as a representation of the decomposition group at $\infty$.
However, the bad factor at $0$ seems complicated to determine for
general $n$. In the case $n=2$, it is easy. Thus, we have the
following complete result for $n=2$, which is the conjectural
Theorem B in [R]. (Our notations are different from those in
Robba).

\bigskip
\noindent {\bf Theorem 0.2.} Suppose $n=2$, $q=p$, and $p$ is an
odd prime. Then $L({\bf G}_m, {\rm Sym}^k({\rm Kl}_2), T)$ is a
polynomial. Its degree is $\frac{k}{2}-[\frac{k}{2p}]$ if $k$ is
even, and $\frac{k+1}{2}-[\frac{k}{2p}+\frac{1}{2}]$ if $k$ is
odd. Moreover, we have the decomposition
$$L({\bf G}_m, {\rm Sym}^k({\rm Kl}_2), T)=P_k(T)
M_k(T)$$ with
\[P_k(T)=\left\{ \begin {array} {ll}
1-T & {\rm if}\; 2\not | k, \\
(1-T)(1-p^{\frac{k}{2}}T)^{m_k}& {\rm if} \; 2|k \;{\rm
and} \; p\equiv 1 \;{\rm mod}\; 4, \\
(1-T)(1+p^{\frac{k}{2}}T)^{n_k}(1-p^{\frac{k}{2}}T)^ {m_k-n_k}&
{\rm if} \; 2|k \;{\rm and} \; p\equiv -1 \;{\rm mod}\; 4,
\end{array}
\right.\] where
\[m_k=\left\{ \begin {array} {ll}
1+[\frac{k}{2p}] & {\rm if}\; 4 | k, \\
{}[\frac{k}{2p}]& {\rm if} \; 4\not |k,
\end{array}
\right.\] and $n_k=[\frac{k}{4p}+\frac{1}{2}]$, and $ M_k$ is a
polynomial satisfying the functional equation
$$M_k(T)=ct^\delta M_k
(\frac{1}{p^{k+1}T}),$$ where $c$ is a nonzero constant (depending
on $k$) and $\delta={\rm deg} M_k$.

\bigskip
Note that the slightly different formula for $n_k$ in the
conjectural Theorem B of [R] is incorrect.

\bigskip
\noindent {\bf Remarks}. Theorem 0.2 only addresses the
$L$-function $L({\bf G}_m, {\rm Sym}^k({\rm Kl}_2), T)$ and the
polynomial $M_k(T)$ from the complex point of view. A more
interesting arithmetic problem is to understand their $p$-adic
properties. For example, the first basic question would be to
determine the $p$-adic Newton polygon of the polynomial $M_k(T)$
with integer coefficients. This is expected to be a difficult
problem. A weak but already non-trivial version is to give an
explicit quadratic lower bound for the $p$-adic Newton polygon of
$M_k(T)$, which is uniform in $k$. Such a uniform quadratic lower
bound is known [W1] in the geometric case of the universal family
of elliptic curves over ${\bf F}_p$, where one considers the
$L$-function of the $k$-th symmetric product of the first relative
$\ell$-adic cohomology which is lisse of rank two outside the
cusps. This latter $L$-function, which is an analogue of the above
$M_k(T)$, is essentially (up to some trivial bad factors) the
Hecke polynomial of the $U_p$-operator acting on the space of
weight $k+2$ cusp forms. From this point of view, we can also ask
for an explicit automorphic interpretation of the polynomial
$M_k(T)$ with integer coefficients. The question on the slope
variation of $M_k(T)$ as $k$ varies $p$-adically is related to the
Gouv\^ea-Mazur type conjecture, see Section $2$ in [W3] for a
simple exposition.

More generally, for any fixed $n$, the sequence of $L$-functions
$L({\bf G}_m, {\rm Sym}^k({\rm Kl}_n), T)$ is $p$-adically
continuous in $k$ as a formal power series in $T$ with $p$-adic
coefficients. In fact, viewed as $p$-adic integers, it is easy to
show that the numbers $\pi_i(\lambda)$ can be re-ordered such that
$\pi_1(\lambda)$ is a $1$-unit and all other $\pi_i(\lambda)$ are
divisible by $p$ (the exact slopes of the $\pi_i(\lambda)$ are
determined by Sperber [Sp]). From this and the Euler product
definition of $L({\bf G}_m, {\rm Sym}^k({\rm Kl}_n), T)$, one
checks that if $k_1=k_2+p^mk_3$ with $k_2$ and $k_3$ non-negative
integers, we have the following congruence:
$$L({\bf G}_m, {\rm Sym}^{k_1}({\rm Kl}_n), T) \equiv ~
L({\bf G}_m, {\rm Sym}^{k_2}({\rm Kl}_n), T) ~({\rm mod}~p^m).$$
For any $p$-adic integer $s\in {\bf Z}_p$, let $k_i$ be an
infinite sequence of strictly increasing positive integers which
converge $p$-adically to $s$. Then, the limit
$$L(n; s,T) =\lim_{i\rightarrow \infty} L({\bf G}_m, {\rm
Sym}^{k_i}({\rm Kl}_n), T)$$ exists as a formal power series in
$T$ with $p$-adic integral coefficients. It is independent of the
choice of the sequence $k_i$ we choose. This power series is
closely related to Dwork's unit root zeta function. It follows
from [W2] that for any $p$-adic integer $s$, the $L$-function
$L(n; s,T)$ is the $L$-function of some infinite rank nuclear
overconvergent $\sigma$-module over ${\bf G}_m$. In particular,
$L(n;s,T)$ is $p$-adic meromorphic in $T$. In fact, $L(n;s,T)$ is
meromorphic in the two variables $(s, T)$ with $|s|_p \leq 1$.
Grosse-Kl\"onne [GK] has extended the $p$-meromorphic continuation
of $L(n; s, T)$ to a larger disk of $s$ with $|s|_p <1 +\epsilon$
for some $\epsilon >0$. Presumably, this two variable $L$-function
$L(n; s, T)$ is related to some type of $p$-adic $L$-functions
over number fields. It would be very interesting to understand the
slopes of the zeros and poles of these $p$-adic meromorphic
$L$-functions. Some explicit partial results were obtained in
[W2], see [W4] for a self-contained exposition of such
$L$-functions in the general case.

\bigskip
This paper is organized as follows. In \S 1, we study the local
monodromy of the Kloosterman sheaf at $\infty$. The main result is
Theorem 1.1 which determines ${\rm Kl}_n$ as a representation of
the decomposition subgroup at $\infty$.  In \S 2, we calculate the
bad factors at $\infty$ of the $L$-functions of the symmetric
products of the Kloosterman sheaf. Using these results, we can
then complete the proof of Theorem 1.1. In \S 3, we calculate the
degrees of the $L$-functions of the symmetric products of the
Kloosterman sheaf. In particular, Theorem 0.1 is proved and some
examples are given.  Finally in \S 4, we study the special case
$n=2$ and prove Theorem 0.2.

\bigskip
\noindent {\bf Acknowledgements.} The research of Lei Fu is
supported by the Qiushi Science \& Technologies Foundation, by the
Fok Ying Tung Education Foundation, by the Transcentury Training
Program Foundation, and by the Project 973. The research of Daqing
Wan is partially supported by the NSF and the NNSF of China
(10128103).

\bigskip
\centerline {\bf 1. Local Monodromy at $\infty$}

\bigskip

In [D1], Deligne constructs a lisse $\overline{\bf Q}_l$-sheaf
${\rm Kl}_n$ on ${\bf G}_m={\bf P}^1-\{0, \infty\}$, which we call
the Kloosterman sheaf. It is lisse of rank $n$, puncturely pure of
weight $n-1$, tamely ramified at $0$, totally wild at $\infty$
with Swan conductor $1$, and for any $x\in {\bf G}_m({\bf
F}_{q^k})={\bf F}_{q^k}^\ast$, we have
$${\rm Tr}(F_x, {\rm Kl}_{n,\bar x})=(-1)^{n-1} {\rm Kl}_n({\bf F}_{q^k},
x),$$ where $F_x$ is the geometric Frobenius elements at the point
$x$. Since ${\rm Kl}_n$ is a lisse sheaf on ${\bf G}_m$, it
corresponds to a galois representation of the function field ${\bf
F}_q(t)$ of ${\bf G}_m$. In this section, we give a detailed study
of ${\rm Kl}_n$ as a representation of the decomposition subgroup
at $\infty$. This result will then be used to study $L$-functions
of symmetric products of ${\rm Kl}_n$.

Before stating the main theorem of this section, let's introduce
some notations. Fix a separable closure $\overline{{\bf F}_q(t)}$
of ${\bf F}_q(t)$. Let $x$ be an element in $\overline {{\bf
F}_q(t)}$ satisfying $x^q-x=t$. Then ${\bf F}_q(t,x)$ is galois
over ${\bf F}_q(t)$. We have a canonical isomorphism
$${\bf F}_q\stackrel {\cong}\to {\rm Gal}({\bf F}_q(t,x)/{\bf
F}_q(t))$$ which sends each $a\in {\bf F}_q$ to the element in
${\rm Gal}({\bf F}_q(t,x)/{\bf F}_q(t))$ defined by $x\mapsto
x+a$. For the additive character $\psi:{\bf F}_q\to \overline {\bf
Q}_l^\ast$, let ${\cal L}_\psi$ be the galois representation
defined by
$${\rm Gal}(\overline{{\bf F}_q(t)}/{\bf F}_q(t))\to
{\rm Gal}({\bf F}_q(t,x)/{\bf F}_q(t))\stackrel {\cong}\to {\bf
F}_q \stackrel {\psi^{-1}}\to \overline {\bf Q}_l^\ast.$$ It is
unramfied outside $\infty$, and totally wild at $\infty$ with Swan
conductor $1$. This galois representation defines a lisse
$\overline {\bf Q}_l$-sheaf on ${\bf A}^1={\bf P}^1-\{\infty\}$
which we still denote by ${\cal L}_\psi$.

Let $\mu_m=\{\mu\in {\bf F}|\mu^m=1\}$ be the subgroup of ${\bf
F}^\ast$ consisting of $m$-th roots of unity. Suppose $m|(q-1)$.
Then $\mu_m$ is contained in ${\bf F}_q$. Let $y$ be an element in
$\overline {{\bf F}_q(t)}$ satisfying $y^m=t$. Then ${\bf
F}_q(t,y)$ is galois over ${\bf F}_q(t)$. We have a canonical
isomorphism
$${\bf \mu}_m \stackrel {\cong}\to {\rm Gal}({\bf F}_q(t,y)/{\bf
F}_q(t))$$ which sends each $\mu \in {\bf \mu}_m $ to the element
in ${\rm Gal}({\bf F}_q(t,y)/{\bf F}_q(t))$ defined by $y\mapsto
\mu y$. For any character $\chi: {\bf \mu}_m\to \overline {\bf
Q}_l^\ast$, let ${\cal L}_\chi$ be the galois representation
defined by
$${\rm Gal}(\overline{{\bf F}_q(t)}/{\bf F}_q(t))\to{\rm Gal}({\bf F}_q(t,y)/{\bf
F}_q(t))\stackrel {\cong}\to {\bf \mu}_m \stackrel {\chi^{-1}}\to
\overline {\bf Q}_l^\ast.$$ It is unramified outside $0$ and
$\infty$, and tamely ramified at $0$ and $\infty$. This galois
representation defines a lisse $\overline {\bf Q}_l$-sheaf on
${\bf G}_m$ which we still denote by ${\cal L}_\chi$.

Let $\theta:{\rm Gal}({\bf F}/{\bf F}_q)\to \overline {\bf
Q}_l^\ast$ be a character of the galois group of the finite field.
Denote by ${\cal L}_\theta$  the galois representation
$${\rm Gal}(\overline {{\bf F}_q(t)}/{\bf F}_q(t))
\to {\rm Gal}({\bf F}/{\bf F}_q) \stackrel {\theta} \to \overline
{\bf Q}_l^\ast.$$  It is unramified everywhere, and hence defines
a lisse $\overline {\bf Q}_l$-sheaf on ${\bf P}^1$ which we still
denote by ${\cal L}_\theta$.

Finally let $\overline {\bf Q}_l\left(\frac {1-n}{2}\right)$ be
the sheaf on ${\rm Spec}{\bf F}_q$ corresponding to the galois
representation of ${\rm Gal}({\bf F}/{\bf F}_q)$ which maps the
geometric Frobenius to $q^{\frac{n-1}{2}}$. For any scheme over
${\bf F}_q$, the inverse image of $\overline {\bf Q}_l\left(\frac
{1-n}{2}\right)$ on this scheme is also denoted by the same
notation.

Now we are ready to state the main theorem of this section.

\bigskip
\noindent {\bf Theorem 1.1.} Suppose $n|(q-1)$. As a
representation of the decomposition subgroup $D_\infty$ at
$\infty$, the Kloosterman sheaf ${\rm Kl}_n$ is isomorphic to
$$[n]_\ast({\cal L}_{\psi_n}\otimes {\cal L}_\chi) \otimes {\cal L}_\theta
\otimes \overline {\bf Q}_l\left(\frac {1-n}{2}\right),$$ where
$[n]:{\bf G}_m\to {\bf G}_m$ is the morphism defined by $x\mapsto
x^n$, $\psi_n$ is the additive character
$$\psi_n(a)=\psi(na),$$ $\chi$ is trivial if $n$ is odd, and $\chi$ is
the (unique) nontrivial character $$\chi:{\bf \mu}_2\to\overline
{\bf Q}_l^\ast$$ if $n$ is even, and $\theta:{\rm Gal}({\bf
F}/{\bf F}_q)\to \overline{\bf Q}_l^\ast$ is a character of ${\rm
Gal}({\bf F}/{\bf F}_q)$ with the following properties:

(1) If $n$ is odd, then $\theta$ is trivial.

(2) If $n$ is even, then $\theta^2$ can be described as follows:
Let $\zeta$ be a primitive $n$-th roots of unity in ${\bf F}_q$.
Fix a square root $\sqrt{\zeta}$ of $\zeta$ in ${\bf F}$. We have
a monomorphism
$${\rm Gal}({\bf F}_q(\sqrt{\zeta})/{\bf F}_q)\hookrightarrow {\bf
\mu}_2$$ defined by sending each $\sigma \in {\rm Gal}({\bf
F}_q(\sqrt{\zeta})/{\bf F}_q)$ to
$\frac{\sigma(\sqrt{\zeta})}{\sqrt{\zeta}}\in\mu_2$. The character
$\theta^2$ is the composition
$${\rm Gal}({\bf F}/{\bf F}_q)\to{\rm Gal}({\bf F}_q(\sqrt{\zeta})/{\bf F}_q)\hookrightarrow
\mu_2\stackrel {\chi^{\frac{n}{2}}}\to \overline {\bf Q}_l^\ast.$$

\bigskip
{\bf Remark}. For even $n$, the above description of $\theta^2$
shows that $\theta^2=1$ if $n\equiv 0$ (mod $4$). If $n \equiv 2$
(mod $4$) and $4|(q-1)$, then since $n|(q-1)$, we must have
$2n|(q-1)$. So $(\sqrt\zeta)^{q-1}=\zeta^{\frac{q-1}{2}}=1$ and
hence $\sqrt{\zeta}\in {\bf F}_q$. So $\theta^2=1$ if $n \equiv 2$
(mod $4$) and $4|(q-1)$. In the remaining case that $n \equiv 2$
(mod $4$) and $4\not | (q-1)$, we have $\sqrt\zeta\not\in {\bf
F}_q$ and $\theta^2$ is a primitive quadratic character. We don't
know how to determine $\theta$ itself completely in the case where
$n$ is even.

Throughout this section, we assume $n|(q-1)$. Then the group
$\mu_n=\{\mu\in {\bf F}|\mu^n=1\}$ of $n$-th roots of unity in
${\bf F}$ is contained in ${\bf F}_q$. Fix an algebraic closure
$\overline {{\bf F}_q(t)}$ of ${\bf F}_q(t)$. The aim of this
section is to determine as much as possible ${\rm Kl}_n$ as a
representation of the decomposition subgroup $D_{\infty}$ at
$\infty$. The proof of Theorem 1.1 will be completed in Section 2.
Before that, we need a series of Lemmas.  Our starting point is
the following result of Katz:

\bigskip
\noindent {\bf Lemma 1.2.} As a representation of the wild inertia
subgroup $P_\infty$ at $\infty$, the Kloosterman sheaf ${\rm
Kl}_n$ is isomorphic to $[n]_\ast{\cal L}_{\psi_n}$.

\bigskip
\noindent {\bf Proof.}  This follows from Propositions 10.1 and
5.6.2 in [K].
\bigskip
\noindent {\bf Lemma 1.3.} Let $G$ be a group, $H$ a subgroup of
$G$ with finite index, and $\rho: H\to {\rm GL}(V)$ a
representation. Suppose $H$ is normal in $G$. For any $g\in G$,
let $\rho_g$ be the composition
$$H\stackrel {{\rm adj}_g}\to H\stackrel \rho\to {\rm GL}(V),$$
where ${\rm adj}_g(h)=g^{-1}hg$. Then the isomorphic class of the
representation $\rho_g$ depends only on the image of $g$ in $G/H$,
and
$${\rm Res}_H {\rm Ind}_H^G \rho\cong \bigoplus_{g\in G/H} \rho_g.
$$

\bigskip
\noindent {\bf Proof.} This is a special case of Proposition 22 on
Page 58 of [S].

\bigskip
\noindent {\bf Lemma 1.4.} We have $$[n]^\ast [n]_\ast {\cal
L}_{\psi_n}\cong \bigoplus_{\mu^n=1} {\cal L}_{\psi_{n\mu}},$$
where $\psi_{n\mu}$ is the additive character
$\psi_{n\mu}(x)=\psi(n\mu x)$.

\bigskip
\noindent {\bf Proof.} Let $y,z$ be elements in $\overline{{\bf
F}_q(t)}$ satisfying $y^n=t$ and $z^q-z=y$. Note that ${\bf
F}_q(z)$ and ${\bf F}_q(y)$ are galois extensions of ${\bf
F}_q(t)$. Let $G={\rm Gal}({\bf F}_q(z)/{\bf F}_q(t))$ and $H={\rm
Gal}({\bf F}_q(z)/{\bf F}_q(y))$. Then $H$ is a normal subgroup of
$G$. Let $\rho$ be the representation
$$H={\rm Gal}({\bf
F}_q(z)/{\bf F}_q(y))\stackrel {\cong} \to {\bf F}_q\stackrel
{\psi_n^{-1}}\to \overline {\bf Q}_l^\ast.$$ Then $[n]^\ast
[n]_\ast {\cal L}_{\psi_n}$ is isomorphic to the composition of
${\rm Res}_H{\rm Ind}_H^G\rho$ with the canonical homomorphism
${\rm Gal}(\overline{{\bf F}_q(t)}/{\bf F}_q(y))\to {\rm Gal}({\bf
F}_q(z)/{\bf F}_q(y))=H$. We have canonical isomorphisms
$$G/H\stackrel {\cong }\to {\rm Gal}({\bf F}_q(y)/{\bf
F}_q(t)) \stackrel {\cong}\to \mu_n.$$ For each $\mu\in \mu_n$,
let $g_\mu\in G={\rm Gal}({\bf F}_q(z)/{\bf F}_q(t))$ be the
element defined by $g_\mu(z)=\mu z$. Then the images of $g_\mu$
($\mu\in \mu_n$) in $G/H$ form a family of representatives of
cosets. By Lemma 1.3, we have $${\rm Res}_H{\rm Ind}_H^G\rho\cong
\bigoplus_{\mu^n=1}\rho_{g_\mu},$$ where $\rho_{g_\mu}$ is the
composition
$$H\stackrel {{\rm adj}_{g_\mu}}\to H \stackrel {\cong }\to {\bf
F}_q \stackrel {\psi_n^{-1}}\to \overline {\bf Q}_l^\ast.$$ One
can verify that we have a commutative diagram
$$\begin{array}{rcl}
H&\stackrel {{\rm adj}_{g_\mu}}\to& H \\
\cong \downarrow &&\downarrow\cong \\
{\bf F}_q &\stackrel{a\mapsto \mu a}\to &{\bf F}_q,
\end{array}$$
where the vertical arrows are the canonical isomorphism $H={\rm
Gal}({\bf F}_q(z)/{\bf F}_q(y))\stackrel {\cong}\to  {\bf F}_q$.
So $\rho_{g_\mu}$ is the composition
$$H \stackrel {\cong }\to {\bf
F}_q \stackrel {\psi_{n\mu}^{-1}}\to \overline {\bf Q}_l^\ast.$$
This proves the lemma.

\bigskip
Before stating the next lemma, let us introduce some notations.
Let $\eta_\infty$ be the generic point of the henselization of
${\bf P}^1$ at $\infty$ and let $\overline\eta_\infty$ be a
geometric point located at $\eta_\infty$. Fix an embedding of
$\overline {{\bf F}_q(t)}$ into the residue field $k(\overline
\eta_\infty)$ of $\overline \eta_\infty$. This defines a
monomorphism ${\rm Gal}(k(\overline
\eta_\infty)/k(\eta_\infty))\hookrightarrow {\rm
Gal}(\overline{{\bf F}_q(t)}/{\bf F}_q(t))$ whose image is the
decomposition subgroup $D_\infty$ of ${\rm Gal}(\overline{{\bf
F}_q(t)}/{\bf F}_q(t))$ at $\infty$. We identify ${\rm
Gal}(k(\overline \eta_\infty)/k(\eta_\infty))$ with $D_\infty$
through this monomorphism. The category of lisse $\overline {\bf
Q}_l$-sheaves on $\eta_\infty$ is equivalent to the category of
$\overline {\bf Q}_l$-representations of $D_\infty$. For
convenience, we denote a lisse sheaf on $\eta_\infty$ and the
corresponding representation of $D_\infty$ by the same symbol. The
morphism $[n]:{\bf G}_m\to {\bf G}_m$ induces a morphism
$\eta_\infty\to \eta_\infty$ which we still denote by $[n]$. 

\bigskip
\noindent {\bf Lemma 1.5.} As a representation of the
decomposition subgroup $D_\infty$ at $\infty$, the Kloosterman
sheaf ${\rm Kl}_n$ is isomorphic to $[n]_\ast({\cal
L}_{\psi_n}\otimes \phi)$ for some tamely ramified one dimensional
representation $\phi:D_\infty \to \overline {\bf Q}_l^\ast$.

\bigskip
\noindent {\bf Proof.} Let $V$ be the stalk of the Kloosterman
sheaf ${\rm Kl}_n$ at $\overline \eta_\infty$. Then $D_\infty$
acts on $V$. Since the Swan conductor of ${\rm Kl}_n$ at $\infty$
is $1$, $V$ is irreducible as a representation of the inertia
group $I_\infty$ and hence irreducible as a representation of
$D_\infty$. Since $n$ is relatively prime to $p$, the morphism
$[n]:{\bf G}_m\to {\bf G}_m$ induces an isomorphism $[n]_\ast:
P_\infty\stackrel\cong\to P_\infty$ on the wild inertia subgroup
$P_\infty$. By Lemma 1.2, $V$ has a one-dimensional subspace $L$,
stable under the action of $P_\infty$, and isomorphic to ${\cal
L}_{\psi_n}$ as a representation of $P_\infty$. By Lemma 1.4, the
restriction of $V$ to $P_\infty$ is not isotypic. Let
$D_\infty^\prime$ be the subgroup of $D_\infty$ consisting of
those elements leaving $L$ stable. Then by Proposition 24 on page
61 of [S], $V$ is isomorphic to ${\rm
Ind}_{D_\infty^\prime}^{D_\infty}(L)$ as a representation of
$D_\infty$. Since the rank of $V$ is $n$, $D_\infty^\prime$ is a
subgroup of $D_\infty$ with index $n$. It defines a finite
extension of degree $n$ over the residue field $k(\eta_\infty)$ of
$\eta_\infty$. Since $D_\infty^\prime$ contains $P_\infty$, this
finite extension is tamely ramified. The degree of inertia of this
finite extension is necessarily $1$. Otherwise, the family of
double cosets $I_\infty\backslash D_\infty/D_\infty^\prime$ would
contain more than one elements and ${\rm Res}_{I_\infty} {\rm
Ind}_{D_\infty^\prime}^{D_\infty}(L)$ would not be irreducible by
Proposition 22 on page 58 of [S]. This contradicts to the fact
that $V$ is irreducible as a representation of $I_\infty$. These
facts imply that the above finite extension is just
$[n]:\eta_\infty\to \eta_\infty$ and $D_\infty^\prime$ is the
image of the monomorphism $[n]_\ast:D_\infty\hookrightarrow
D_\infty$ induced by the morphism $[n]$. Since $L$ is isomorphic
${\cal L}_{\psi_n}$ as a representation of $P_\infty$, $L$ is
isomorphic to ${\cal L}_{\psi_n}\otimes L^\prime$ as a
representation of $D_\infty^\prime$, where $L^\prime$ is a
representation of $D_\infty^\prime$ of rank $1$ which is trivial
when restricted to $P_\infty$, that is, $L^\prime$ is tamely
ramified. Composing with the isomorphism
$[n]_\ast:D_\infty\stackrel {\cong}\to D_\infty^\prime$,
$L^\prime$ defines a tamely ramified one-dimensional
representation $\phi:D_\infty\to \overline {\bf Q}_l^\ast $, and
$V$ is isomorphic to $[n]_\ast({\cal L}_{\psi_n}\otimes \phi)$.

\bigskip
\noindent {\bf Lemma 1.6.} Keep the notation of Lemma 1.5. As a
representation of the inertia subgroup $I_\infty$, $\phi$ is
isomorphic to ${\cal L}_\chi$, where $\chi$ is trivial if $n$ is
odd, and $\chi$ is the unique nontrivial character $\chi:\mu_2\to
\overline {\bf Q}_l^\ast$ if $n$ is even.

\bigskip
\noindent {\bf Proof.} First recall that making the base extension
from ${\bf F}_q$ to ${\bf F}$ has no effect on the inertia
subgroup $I_\infty$. So we can work over the base ${\bf F}$. By
the Appendix of [ST], the representation $\phi:D_\infty\to
\overline{\bf Q}_l^\ast$ is quasi-unipotent when restricted to
$I_\infty$, and hence has finite order when restricted to
$I_\infty$ (since the representation is one-dimensional). Since
$\phi$ is tamely ramified, there exists a positive integer $m$
relatively prime to $p$ so that as a representation of $I_\infty$,
$\phi$ is isomorphic to the restriction to $I_\infty$ of the
galois representation
$${\rm Gal}(
{{\bf F}(\sqrt[m]{y})}/{\bf F}(y)) \stackrel \cong\to \mu_m
\stackrel {\chi^{-1}} \to \overline {\bf Q}_l^\ast$$ for some
primitive character $\chi: \mu_m\to \overline {\bf Q}_l^\ast$.

Let us calculate ${\rm det}([n]_\ast ({\cal L}_{\psi_n}\otimes
{\cal L}_{\chi}))$. Let $y$, $z$, $w$ be elements in $\overline
{{\bf F}(t)}$ satisfying $y^n=t$, $z^q-z=y$, and $w^m=y$. Then
${\bf F}(z, w)$ and ${\bf F}(y)$ are galois extensions of ${\bf
F}(t)$. Let $G={\rm Gal}({\bf F}(z,w)/{\bf F}(t))$ and $H={\rm
Gal}({\bf F}(z,w)/{\bf F}(y))$. Then $H$ is normal in $G$, and we
have canonical isomorphisms
$$G/H\stackrel\cong\to {\rm Gal}({\bf F}(y)/{\bf F}(t))\stackrel
\cong \to \mu_n.$$ We have an isomorphism
$${\bf F}_q \times \mu_m \stackrel \cong\to
H={\rm Gal}({\bf F}(z,w)/{\bf F}(y))$$ which maps $(a,\mu)\in {\bf
F}_q\times \mu_m$ to the element $g_{(a,\mu)}\in {\rm Gal}({\bf
F}(z,w)/{\bf F}(y))$ defined by $g_{(a,\mu)}(z)=z+a$ and
$g_{(a,\mu)}(w)=\mu w$. Let $\omega:H\to \overline {\bf Q}_l^\ast$
be the character defined by
$$\omega(g_{(a,\mu)})=\psi_n(-a)\chi(\mu^{-1}).$$
Then $[n]_\ast ({\cal L}_{\psi_n}\otimes {\cal L}_{\chi})$ is just
the composition of ${\rm Ind}_H^G(\omega)$ with the canonical
homomorphism ${\rm Gal}(\overline {{\bf F}(t)}/{\bf F}(t))\to {\rm
Gal}({\bf F}(z,w)/{\bf F}(t))=G$. Let $\zeta$ be a primitive
$n$-th root of unity in ${\bf F}$. Choose an $m$-th root
$\sqrt[m]{\zeta}$ of $\zeta$. (Then $\sqrt[m]\zeta$ is a primitive
$mn$-th root of unity). Let $g$ be the element in $G={\rm
Gal}({\bf F}(z,w)/{\bf F}(t))$ defined by $g(z)=\zeta z$ and
$g(w)=\sqrt[m]{\zeta} w$. Then the image of $g$ in $G/H$ is a
generator of the cyclic group $G/H$. So $G$ is generated by
$g_{(a,\mu)}\in H$ $((a,\mu)\in {\bf F}_q\times \mu_m)$ and $g$.
By Lemma 1.3, we have
$${\rm Res}_H{\rm Ind}_H^G (\omega)=\bigoplus_{i=0}^{n-1}
\omega_{g^i},$$ where $\omega_{g^i}$ is the composition
$$H\stackrel {{\rm adj}_{g^i}}\to H \stackrel \omega\to \overline{\bf Q}_l^\ast.$$
One can verify that
$$\omega_{g^i}(g_{(a,\mu)})=\psi_n(-\zeta^i a)\chi(\mu^{-1}).$$
So
\begin{eqnarray*}
({\rm det}({\rm Ind}_H^G
(\omega)))(g_{(a,\mu)})&=&\prod_{i=0}^{n-1}
(\psi_n(-\zeta^i a)\chi(\mu^{-1}))\\
&=&\psi_n(-(\sum_{i=0}^{n-1}\zeta^i)a)\chi(\mu^{-n})\\
&=& \psi_n(0)\chi(\mu^{-n})\\
&=&\chi(\mu^{-n})
\end{eqnarray*}
We have
\begin{eqnarray*}
({\rm det})({\rm Ind}_H^G (\omega)))(g) &=& {\rm det}\left(
\begin{array}{cccc}
&1&& \\
&&\ddots& \\
&&&1 \\
\omega(g^n)&&&
\end{array}
\right)\\
&=&(-1)^{n+1}\omega(g^n).\\
\end{eqnarray*}
One can verify $g^n=g_{(0, (\sqrt[m]{\zeta})^{n})}$. So
$\omega(g^n)=\chi((\sqrt[m]{\zeta})^{-n})$ and hence
$$({\rm det})({\rm Ind}_H^G
(\omega)))(g)=(-1)^{n+1}\chi((\sqrt[m]{\zeta})^{-n}).$$ By Lemma
1.5, ${\rm det}({\rm Kl}_n)$ is isomorphic to ${\rm
det}([n]_\ast({\cal L}_{\psi_n}\otimes\phi))$ as a representation
of $I_\infty$. By [K] 7.4.3, ${\rm det}({\rm Kl}_n)$ is
geometrically constant. On the other hand, as representations of
$I_\infty$, we have ${\rm det}([n]_\ast({\cal
L}_{\psi_n}\otimes\phi))\cong {\rm det}([n]_\ast({\cal
L}_{\psi_n}\otimes{\cal L}_{\chi}))$, and ${\rm
det}([n]_\ast({\cal L}_{\psi_n}\otimes{\cal L}_{\chi}))$ is
isomorphic to the composition of ${\rm det}({\rm
Ind}_H^G(\omega))$ with the canonical homomorphism ${\rm
Gal}(\overline {{\bf F}(t)}/{\bf F}(t))\to {\rm Gal}({\bf
F}(z,w)/{\bf F}(t))=G$. So ${\rm det}({\rm Ind}_H^G(\omega))$ is
trivial as a representation of $I_\infty$. Hence
$\chi(\mu^{-n})=1$ for all $\mu\in \mu_m$ and
$\chi((\sqrt[m]{\zeta})^{-n})=(-1)^{n+1}$. Since $\zeta$ is a
primitive $n$-th root of unity, $(\sqrt[m]{\zeta})^{-n}$ is a
primitive $m$-th root of unity. This implies that the order $m$ of
$\chi$ is at most $2$. If $n$ is odd, the relation
$\chi((\sqrt[m]{\zeta})^{-n})=(-1)^{n+1}=1$ implies that $\chi$ is
trivial. If $n$ is even, the relation
$\chi((\sqrt[m]{\zeta})^{-n})=(-1)^{n+1}=-1$ implies that the
order $m$ of $\chi$ is exactly $2$. This finishes the proof of the
lemma.

\bigskip
\noindent {\bf Lemma 1.7.} As a representation of $D_\infty$, the
Kloosterman sheaf ${\rm Kl}_n$ is isomorphic to
$$[n]_\ast({\cal L}_{\psi_n}\otimes {\cal L}_\chi)\otimes {\cal
L}_\theta\otimes \overline {\bf Q}_l\left(\frac{1-n}{2}\right),$$
where $\chi$ is trivial if $n$ is odd, and $\chi$ is the
nontrivial character $\chi:\mu_2\to \overline {\bf Q}_l^\ast$ if
$n$ is even, and $\theta:{\rm Gal}({\bf F}/{\bf F}_q)\to
\overline{\bf Q}_l^\ast$ is a character of ${\rm Gal}({\bf F}/{\bf
F}_q)$.

\bigskip
\noindent{\bf Proof.} By Lemma 1.5, as a representation of
$D_\infty$, ${\rm Kl}_n$ is isomorphic to $[n]_\ast({\cal
L}_{\psi_n}\otimes\phi)$ for some character
$\phi:D_\infty\to\overline {\bf Q}_l^\ast$. By Lemma 1.6, $\phi$
is isomorphic to ${\cal L}_\chi$ when restricted to $I_\infty$. So
as a representation of $D_\infty$, $\phi$ is isomorphic to ${\cal
L}_\chi\otimes {\cal L}_\theta\otimes {\bf
Q}_l\left(\frac{1-n}{2}\right)$ for some character $\theta:{\rm
Gal}({\bf F}/{\bf F}_q)\to \overline{\bf Q}_l^\ast$. So as a
representation of $D_\infty$, ${\rm Kl}_n$ is isomorphic to
$$[n]_\ast\left({\cal L}_{\psi_n}\otimes {\cal L}_\chi\otimes
{\cal L}_\theta\otimes \overline {\bf
Q}_l\left(\frac{1-n}{2}\right)\right)\cong [n]_\ast({\cal
L}_{\psi_n}\otimes {\cal L}_\chi)\otimes {\cal L}_\theta\otimes
\overline {\bf Q}_l\left(\frac{1-n}{2}\right).$$ Lemma 1.7 is
proved.

\bigskip
\noindent{\bf Lemma 1.8.} Keep the notation in Lemma 1.7. The
character $\theta:{\rm Gal}({\bf F}/{\bf F}_q)\to \overline{\bf
Q}_l^\ast$ has the following properties:

(1) If $n$ is odd, then $\theta^n$ is trivial.

(2) If $n$ is even, then $\theta^2$ can be described as follows:
Let $\zeta$ be a primitive $n$-th roots of unity in ${\bf F}_q$.
Fix a square root $\sqrt{\zeta}$ of $\zeta$ in ${\bf F}$. We have
a monomorphism
$${\rm Gal}({\bf F}_q(\sqrt{\zeta})/{\bf F}_q)\hookrightarrow {\bf
\mu}_2$$ defined by sending each $\sigma \in {\rm Gal}({\bf
F}_q(\sqrt{\zeta})/{\bf F}_q)$ to
$\frac{\sigma(\sqrt{\zeta})}{\sqrt{\zeta}}\in\mu_2$. The character
$\theta^2$ is the composition
$${\rm Gal}({\bf F}/{\bf F}_q)\to {\rm Gal}({\bf F}_q(\sqrt{\zeta})/{\bf F}_q)\hookrightarrow
\mu_2\stackrel {\chi^{\frac{n}{2}}}\to \overline {\bf Q}_l^\ast,$$
where $\chi:\mu_2\to \overline{\bf Q}_l^\ast$ is the nontrivial
character on $\mu_2$.

(3) If $n$ is odd and $p=2$, then $\theta$ is trivial.

\bigskip
\noindent {\bf Proof.} Suppose $n$ is odd. By Lemma 1.7, ${\rm
Kl}_n$ is isomorphic to $[n]_\ast {\cal L}_{\psi_n}\otimes {\cal
L}_\theta\otimes \overline {\bf Q}_l\left(\frac{1-n}{2}\right)$ as
a representation of $D_\infty$. Using the same method as in the
proof of Lemma 1.6 (but working over the base ${\bf F}_q$), one
can show that ${\rm det}([n]_\ast{\cal L}_{\psi_n})$ is trivial.
So $${\rm det}\left([n]_\ast {\cal L}_{\psi_n}\otimes {\cal
L}_\theta\otimes \overline {\bf
Q}_l\left(\frac{1-n}{2}\right)\right)={\cal
L}_{\theta^n}\otimes\overline{\bf Q}_l
\left(\frac{n(1-n)}{2}\right).$$ By [K] 7.4.3, we have ${\rm
det}({\rm Kl}_n)=\overline{\bf Q}_l\left(\frac{n(1-n)}{2}\right)$.
So $\theta^n$ is trivial. This proves part (1) of Lemma 1.8.

Suppose $n$ is even. In this case,  ${\rm Kl}_n$ is isomorphic to
$[n]_\ast({\cal L}_{\psi_n}\otimes {\cal L}_\chi)\otimes {\cal
L}_\theta\otimes \overline {\bf Q}_l\left(\frac{1-n}{2}\right)$ as
a representation of $D_\infty$, where $\chi$ is the nontrivial
character $\chi:\mu_2\to \overline {\bf Q}_l^\ast$. By [K] 4.1.11
and 4.2.1, there exists a perfect skew-symmetric pairing ${\rm
Kl}_n\times {\rm Kl}_n\to \overline {\bf Q}_l(1-n)$, and any such
pairing invariant under the action of $I_\infty$ coincides with
this one up to a scalar. So the $I_\infty$-coinvariant space
$(\bigwedge^2 {\rm Kl}_n)_{I_\infty}$ of $\bigwedge^2 {\rm Kl}_n$
is isomorphic to $\overline {\bf Q}_l(1-n)$ as a representation of
${\rm Gal}({\bf F}/{\bf F}_q)$. Hence $(\bigwedge^2[n]_\ast({\cal
L}_{\psi_n}\otimes{\cal L}_\chi))_{I_\infty}\otimes {\cal
L}_{\theta^2}\otimes \overline {\bf Q}_l(1-n)$ is isomorphic to
$\overline {\bf Q}_l(1-n)$. Note that $\bigwedge^2 [n]_\ast({\cal
L}_{\psi_n}\otimes{\cal L}_\chi)$ is a semisimple representation
of $I_\infty$ since it has finite monodromy. So the canonical
homomorphism
$$(\bigwedge^2 [n]_\ast({\cal L}_{\psi_n}\otimes{\cal
L}_\chi))^{I_\infty}\to(\bigwedge^2 [n]_\ast({\cal
L}_{\psi_n}\otimes{\cal L}_\chi))_{I_\infty}$$ is an isomorphism.
We will show that as a representation of ${\rm Gal}({\bf F}/{\bf
F}_q)$, $(\bigwedge^2 [n]_\ast({\cal L}_{\psi_n}\otimes{\cal
L}_\chi))^{I_\infty}$ is isomorphic to the composition
$${\rm Gal}({\bf F}/{\bf F}_q)\to{\rm Gal}({\bf F}_q(\sqrt{\zeta})/{\bf F}_q)\hookrightarrow
\mu_2\stackrel {\chi^{\frac{n}{2}}}\to \overline {\bf Q}_l^\ast.$$
This will prove part (2) of Lemma 1.8.

Let $\zeta$ be a primitive $n$-th root of unity in ${\bf F}_q$.
Fix a square root $\sqrt\zeta$ of $\zeta$ in ${\bf F}$. Let $y$,
$z$, $w$ be elements in $\overline {{\bf F}_q(t)}$ satisfying
$y^n=t$, $z^q-z=y$, and $w^2=y$. Then ${\bf F}_q(z, w,\sqrt\zeta)$
and ${\bf F}_q(y)$ are galois extensions of ${\bf F}_q(t)$. Let
$G={\rm Gal}({\bf F}_q(z,w,\sqrt\zeta)/{\bf F}_q(t))$ and $H={\rm
Gal}({\bf F}_q(z,w,\sqrt\zeta)/{\bf F}_q(y))$. Then $H$ is normal
in $G$, and we have canonical isomorphisms
$$G/H\stackrel\cong \to {\rm Gal}({\bf F}_q(y)/{\bf F}_q(t))\stackrel
\cong \to \mu_n.$$ Consider the case where $\sqrt\zeta$ does not
lie in ${\bf F}_q$. We then have an isomorphism
$${\bf F}_q \times \mu_2\times \mu_2 \stackrel \cong\to
H={\rm Gal}({\bf F}_q(z,w,\sqrt\zeta)/{\bf F}_q(y))$$ which maps
$(a,\mu',\mu'')\in {\bf F}_q\times \mu_2\times \mu_2$ to the
element $g_{(a,\mu',\mu'')}\in {\rm Gal}({\bf
F}_q(z,w,\sqrt\zeta)/{\bf F}_q(y))$ defined by
$g_{(a,\mu',\mu'')}(z)=z+a$, $g_{(a,\mu',\mu'')}(w)=\mu' w$, and
$g_{(a,\mu',\mu'')}(\sqrt\zeta)=\mu''\sqrt \zeta$. (In the case
where $\sqrt\zeta$ lies in ${\bf F}_q$, we have ${\bf F}_q(z,
w,\sqrt\zeta)={\bf F}_q(z,w)$, and we have an isomorphism
$${\bf F}_q \times \mu_2\stackrel \cong\to
H={\rm Gal}({\bf F}_q(z,w,\sqrt\zeta)/{\bf F}_q(y))$$ which maps
$(a,\mu')\in {\bf F}_q\times \mu_2$ to the element
$g_{(a,\mu')}\in {\rm Gal}({\bf F}_q(z,w,\sqrt\zeta)/{\bf
F}_q(y))$ defined by $g_{(a,\mu')}(z)=z+a$ and
$g_{(a,\mu')}(w)=\mu' w$. All the following argument works for
this case with slight modification. We leave to the reader to
treat this case.) Let $\omega:H\to \overline {\bf Q}_l^\ast$ be
the character defined by
$$\omega(g_{(a,\mu',\mu'')})=\psi_n(-a)\chi(\mu'^{-1}).$$
Then $[n]_\ast ({\cal L}_{\psi_n}\otimes {\cal L}_{\chi})$ is just
the composition of ${\rm Ind}_H^G(\omega)$ with the canonical
homomorphism ${\rm Gal}(\overline {{\bf F}_q(t)}/{\bf F}_q(t))\to
{\rm Gal}({\bf F}_q(z,w,\sqrt\zeta)/{\bf F}_q(t))=G$. Let $g$ be
the element in $G={\rm Gal}({\bf F}_q(z,w,\sqrt\zeta)/{\bf
F}_q(t))$ defined by $g(z)=\zeta z$, $g(w)=\sqrt{\zeta} w$ and
$g(\sqrt\zeta)=\sqrt\zeta$. Then the image of $g$ in $G/H$ is a
generator of the cyclic group $G/H$. So $G$ is generated by
$g_{(a,\mu',\mu'')}\in H$ $((a,\mu',\mu'')\in {\bf F}_q\times
\mu_2\times\mu_2)$ and $g$. By Lemma 1.3, we have
$${\rm Res}_H{\rm Ind}_H^G (\omega)=\bigoplus_{i=0}^{n-1}
\omega_{g^i},$$ where $\omega_{g^i}$ is the composition
$$H\stackrel {{\rm adj}_{g^i}}\to H \stackrel \omega\to \overline{\bf Q}_l^\ast.$$
One can verify that
$$\omega_{g^i}(g_{(a,\mu',\mu'')})=\psi_n(-\zeta^i a)\chi(\mu'^{-1}\mu''^{-i}).$$
Let $V$ be the stalk of $[n]_\ast({\cal L}_{\psi_n}\otimes{\cal
L}_\chi)$ at the geometric point $\overline \eta_\infty$. The
above calculation shows that there exists a basis $\{e_0,\ldots,
e_{n-1}\}$ of $V$ such that $ge_0=e_1, \; ge_1=e_2,\ldots,
ge_{n-2}=e_{n-1}$, and for any $g_{(a,\mu',\mu'')}\in H$
($(a,\mu',\mu'')\in {\bf F}_q\times \mu_2\times \mu_2$), we have
$$g_{(a,\mu',\mu'')}
(e_i)=\omega_{g^i}(g_{(a,\mu',\mu'')})e_i=\psi_n(-\zeta^i
a)\chi(\mu'^{-1}\mu''^{-i})e_i.$$ One can verify
$g^n=g_{(0,(\sqrt\zeta)^n,1)}$. So
$$g^ne_0=\chi((\sqrt\zeta)^{-n})e_0=-e_0.$$
Consider the element $\sum\limits_{i=0}^{\frac{n}{2}-1} e_i\wedge
e_{i+\frac{n}{2}}\in \bigwedge^2 V$. We have
\begin{eqnarray*}
g(e_0\wedge e_{\frac{n}{2}}+\cdots +e_{\frac{n}{2}-1}\wedge
e_{n-1})&=&e_1\wedge e_{1+\frac{n}{2}}+\cdots +
e_{\frac{n}{2}-1}\wedge
e_{n-1} + e_{\frac{n}{2}}\wedge g^ne_0\\
&=& e_0\wedge e_{\frac{n}{2}}+\cdots +e_{\frac{n}{2}-1}\wedge
e_{n-1}
\end{eqnarray*}
and
\begin{eqnarray*}
&&g_{(a,\mu',\mu'')}(\sum\limits_{i=0}^{\frac{n}{2}-1} e_i\wedge
e_{i+\frac{n}{2}})\\
&=& \sum\limits _{i=0}^{\frac{n}{2}-1} \psi_n(-\zeta^i
a)\chi(\mu'^{-1}\mu''^{-i}) \psi_n(-\zeta^{i+\frac{n}{2}}
a)\chi(\mu'^{-1}\mu''^{-(i+\frac{n}{2})})e_i\wedge e_{i+\frac{n}{2}}\\
&=& \sum\limits _{i=0}^{\frac{n}{2}-1}
\psi_n(-\zeta^ia(1+\zeta^{\frac{n}{2}}))\chi(\mu'^{-2}\mu''^{-2i-\frac{n}{2}})
e_i\wedge e_{i+\frac{n}{2}}\\
&=&\chi(\mu'')^{\frac{n}{2}}\sum\limits_{i=0}^{\frac{n}{2}-1}
e_i\wedge e_{i+\frac{n}{2}}
\end{eqnarray*}
since $\zeta^{\frac{n}{2}}=-1$ and $\chi$ is of order $2$. In
particular, $g$ and $g_{(a,\mu',1)}$ act trivially on
$\sum\limits_{i=0}^{\frac{n}{2}-1} e_i\wedge e_{i+\frac{n}{2}}$.
So $\sum\limits_{i=0}^{\frac{n}{2}-1} e_i\wedge e_{i+\frac{n}{2}}$
lies in $(\bigwedge^2 V)^{I_\infty}$. Note that
$\sum\limits_{i=0}^{\frac{n}{2}-1} e_i\wedge e_{i+\frac{n}{2}}$
spans $(\bigwedge^2 V)^{I_\infty}$ as the latter space is one
dimensional. On the other hand, for any $\mu''\in \mu_2$, we have
$g_{(0,1,\mu'')}(\sum\limits_{i=0}^{\frac{n}{2}-1} e_i\wedge
e_{i+\frac{n}{2}})=\chi(\mu'')^{\frac{n}{2}}\sum\limits_{i=0}^{\frac{n}{2}-1}
e_i\wedge e_{i+\frac{n}{2}}$. So as a representation ${\rm
Gal}({\bf F}/{\bf F}_q)$, $(\bigwedge^2V)^{I_\infty}$ is
isomorphic to the composition
$${\rm Gal}({\bf F}/{\bf F}_q)\to{\rm Gal}({\bf F}_q(\sqrt{\zeta})/{\bf F}_q)\hookrightarrow
\mu_2\stackrel {\chi^{\frac{n}{2}}}\to \overline {\bf Q}_l^\ast.$$
This finishes the proof of part (2) of Lemma 1.8.

Finally suppose $n$ is odd and $p=2$. Then ${\rm Kl}_n$ is
isomorphic to $[n]_\ast {\cal L}_{\psi_n}\otimes {\cal
L}_\theta\otimes \overline {\bf Q}_l\left(\frac{1-n}{2}\right)$ as
a representation of $D_\infty$. By [K] 4.1.11 and 4.2.1, there
exists a perfect symmetric pairing ${\rm Kl}_n\times {\rm Kl}_n\to
\overline {\bf Q}_l(1-n)$, and any such pairing invariant under
the action of $I_\infty$ coincides with this one up to a scalar.
So the $I_\infty$-coinvariant space $({\rm Sym}^2 ({\rm
Kl}_n))_{I_\infty}$ of ${\rm Sym}^2 ({\rm Kl}_n)$ is isomorphic to
$\overline {\bf Q}_l(1-n)$ as a representation of ${\rm Gal}({\bf
F}/{\bf F}_q)$. Hence $({\rm Sym}^2([n]_\ast{\cal
L}_{\psi_n}))_{I_\infty}\otimes {\cal L}_{\theta^2}\otimes
\overline {\bf Q}_l(1-n)$ is isomorphic to $\overline {\bf
Q}_l(1-n)$. Note that ${\rm Sym}^2 ([n]_\ast {\cal L}_{\psi_n})$
is a semisimple representation of $I_\infty$ since it has finite
monodromy. So the canonical homomorphism
$$({\rm Sym}^2 ([n]_\ast {\cal L}_{\psi_n}))^{I_\infty}\to
({\rm Sym}^2 ([n]_\ast{\cal L}_{\psi_n}))_{I_\infty}$$ is an
isomorphism. One can show that as a representation of ${\rm
Gal}({\bf F}/{\bf F}_q)$, $({\rm Sym}^2 ([n]_\ast{\cal
L}_{\psi_n}))^{I_\infty}$ is trivial. (Using the notation in the
previous paragraph, one verifies $e_0^2+e_1^2+\cdots +e_{n-1}^2$
is a generator of $({\rm Sym}^2V)^{I_\infty}$ and the geometric
Frobenius acts trivially on this vector.) So $\theta^2$ is
trivial. By (1), $\theta^n$ is also trivial. So $\theta$ must be
trivial. This finishes the proof of part (3) of Lemma 1.8. The
proof of Lemma 1.8 is complete.

\bigskip
\bigskip
\centerline {\bf 2. The bad factors of the $L$-functions of
symmetric products}

\bigskip
\bigskip
Let $\zeta$ be a primitive $n$-th root of unity in ${\bf F}$. For
each positive integer $k$, let $S_k(n,p)$ be the set of $n$-tuples
$(j_0, \cdots, j_{n-1})$ of non-negative integers satisfying
$j_0+j_1+\cdots +j_{n-1}=k$ and $j_0+j_1\zeta +\cdots
+j_{n-1}\zeta^{n-1}=0$ in ${\bf F}$. Let $\sigma$ denote the
cyclic shifting operator
$$\sigma(j_0, \cdots, j_{n-1}) = (j_{n-1}, j_0, \cdots, j_{n-2}).$$
It is clear that the set $S_k(n,p)$ is $\sigma$-stable.

Let $V$ be a $\overline{\bf Q}_l$-vector space of dimension $n$
with basis $\{e_0, \cdots, e_{n-1}\}$. For an $n$-tuple
$j=(j_0,\cdots, j_{n-1})$ of non-negative integers such that
$j_0+\cdots +j_{n-1}=k$, write
$$e^j =e_0^{j_0}e_1^{j_1}\cdots e_{n-1}^{j_{n-1}}$$
as an element of ${\rm Sym}^kV$. For such an $n$-tuple $j$, we
define
$$v_j = \sum_{i=0}^{n-1} (-1)^{j_{n-1}+\cdots +j_{n-i}} e^{\sigma^i(j)}.$$
This is an element of ${\rm Sym}^kV$. If $k=j_0+j_1+\cdots+
j_{n-1}$ is even, then we have $v_{\sigma(j)}=(-1)^{j_{n-1}}v_j$,
and hence the subspace spanned $v_j$ depends only on the
$\sigma$-orbit of $j$. Let $a_k(n,p)$ be the number of
$\sigma$-orbits in $S_k(n,p)$. When $k$ is even, let $b_k(n,p)$ be
the number of those $\sigma$-orbits in $S_k(n,p)$ such that the
subspace spanned by the orbit is not zero.

\bigskip
\noindent{\bf Lemma 2.1.} Suppose $(n,p)=1$. We have
\[{\rm dim}({\rm Sym}^k ({\rm Kl})_n)^{I_\infty}=\left\{
\begin{array}{ll}
a_k(n,p)& \hbox{ if } n \hbox{ is odd,}\\
0&\hbox { if } n \hbox { is even and } k \hbox { is odd,}\\
b_k(n,p)& \hbox { if } n \hbox { and } k \hbox { are both even.}
\end{array}
\right.\]

\bigskip
\noindent {\bf Proof.} Since $(n,p)=1$ and the inertia subgroup
$I_\infty$ does not change if we make base change from ${\bf F}_q$
to its finite extensions, we may assume that $n|(q-1)$.

We use the notations in the proof of Lemma 1.6. Recall that we
have $m=1$ if $n$ is odd, and $m=2$ if $n$ is even. Let $\beta_m$
be a primitive $m$-th root of unity in ${\bf Q}$. Thus $\beta_m =
1$ if $n$ is odd, and $\beta_m=-1$ if $n$ is even. Let $V$ be the
stalk of $[n]_\ast({\cal L}_{\psi_n}\otimes{\cal L}_\chi)$ at the
geometric point $\overline \eta_\infty$. There exists a basis
$\{e_0,\ldots, e_{n-1}\}$ of $V$ such that $ge_0=e_1, \;
ge_1=e_2,\ldots, ge_{n-2}=e_{n-1}, ge_{n-1} =\beta_m e_0$, and for
any $g_{(a,\mu)}\in H$ ($(a,\mu)\in {\bf F}_q\times \mu_m$), we
have
$$g_{(a,\mu)}
(e_i)=\omega_{g^i}(g_{(a,\mu)})e_i=\psi_n(-\zeta^i a)
\chi(\mu^{-1}) e_i.$$ A basis for ${\rm Sym}^kV$ is
$\{e_0^{j_0}e_1^{j_1}\cdots e_{n-1}^{j_{n-1}}\}$, where $j_i$ are
non-negative integers satisfying $|j|=\sum\limits_{i=0}^{n-1} j_i
=k$. Suppose $v=\sum\limits_{|j|=k} a_je_0^{j_0}e_1^{j_1}\cdots
e_{n-1}^{j_{n-1}}$ lies in $({\rm Sym}^kV)^{I_\infty}$. Then $g$
and $g_{(a,1)}$ $(a\in {\bf F}_q)$ act trivially on it. We have
$$g(v) = g(\sum\limits_{|j|=k} a_je_0^{j_0}e_1^{j_1}\cdots e_{n-1}^{j_{n-1}})
=\sum\limits_{|j|=k}
\beta_m^{j_{n-1}}a_je_0^{j_{n-1}}e_1^{j_0}\cdots
e_{n-1}^{j_{n-2}}.$$ From $g(v)=v$, we get
$$a_j =\beta_m^{j_{n-1}}a_{\sigma(j)} =
\beta_m^{j_{n-1}+j_{n-2}}a_{\sigma^2(j)} =\cdots =
\beta_m^{j_{n-1}+j_{n-2}+\cdots+j_0}a_{\sigma^n(j)}=\beta_m^k
a_j.$$ If $n$ is even and $k$ is odd, then $\beta_m^k=-1$. The
above relation then shows that $a_j=0$. So $({\rm
Sym}^kV)^{I_\infty}=0$ in this case.

Now assume that either $n$ is odd or $k$ is even. The above vector
$v$ is a linear combination of the following vectors
$$v_j = \sum_{i=0}^{n-1} \beta_m^{j_{n-1}+\cdots +j_{n-i}} e^{\sigma^i(j)}.$$
The $g$-invariant subspace $({\rm Sym}^kV)^{g}$ is thus spanned by
the vectors $v_j$ with $|j|=k$, where $j$ runs only over
$\sigma$-orbits. On the other hand, one computes that
\begin{eqnarray*}
g_{(a,1)}(v) &=& g_{(a,1)}(\sum_{|j|=k}
a_je_0^{j_0}e_1^{j_1}\cdots
e_{n-1}^{j_{n-1}})\\
& =&\sum_{|j|=k} \psi_n(-a(j_0+j_1\zeta+\cdots +
j_{n-1}\zeta^{n-1})) a_j e_0^{j_0}e_1^{j_1}\cdots
e_{n-1}^{j_{n-1}}
\end{eqnarray*}
for all $a\in {\bf F}_q$. Since $g_{(a,1)}(v)=v$, if $a_j$ is
non-zero, then we must have $j_0+j_1\zeta+\cdots +
j_{n-1}\zeta^{n-1} =0$ in ${\bf F}$, that is, $j\in S_k(n,p)$.
Thus, we have proved that the inertia invariant $({\rm
Sym}^kV)^{I_\infty}$ is spanned by the vectors $v_j$ where $j$
runs over the $\sigma$-orbits of $S_k(n,p)$.

If $n$ is odd, then $\beta_m=1$ and each $v_j$ is non-zero. As $j$
runs over the $\sigma$-orbits of $S_k(n,p)$, the vectors $v_j$ are
clearly independent. We thus have ${\rm dim}({\rm
Sym}^kV)^{I_\infty}=a_k(n,p)$.

If both $n$ and $k$ are even, then $\beta_m=-1$. In this case,
some of the vectors $v_j$ can be zero. The remaining non-zeros
vectors $v_j$, as $j$ runs over the $\sigma$-orbits of $S_k(n,p)$,
will be linearly independent. We thus have ${\rm dim}({\rm
Sym}^kV)^{I_\infty}=b_k(n,p)$. The lemma is proved.

\bigskip
\noindent {\bf Lemma 2.2.} Suppose $p$ is odd and either $n$ is
odd or $n=2$. Let $j:{\bf G}_m\hookrightarrow {\bf P}^1$ be the
open immersion. Then the $L$-functions $L({\bf G}_m, {\rm
Sym}^k({\rm Kl}_n), T)$ and $L({\bf P}^1, j_\ast({\rm Sym}^k({\rm
Kl}_n)), T)$ are polynomials. All the reciprocal roots of the
polynomial $L({\bf P}^1, j_\ast({\rm Sym}^k({\rm Kl}_n)), T)$ have
weight $k(n-1)+1$. The polynomial $L({\bf G}_m, {\rm Sym}^k({\rm
Kl}_n), T)$ has coefficients in ${\bf Z}$.

\bigskip
\noindent{\bf Proof.} By Grothendieck's formula for $L$-functions,
we have
$$L({\bf
P}^1, j_\ast({\rm Sym}^k({\rm Kl}_n)), T)=\prod_{i=0}^2 {\rm
det}(1-FT, H^i({\bf P}^1\otimes {\bf F},j_\ast({\rm Sym}^k({\rm
Kl}_n)))^{(-1)^{i+1}}.$$ Under our conditions on $p$ and $n$, the
global monodromy group of ${\rm Kl}_n$ is ${\rm SL}(n)$ by [K]
11.1. In particular ${\rm Sym}^k({\rm Kl}_n)$ is irreducible as a
representation of the geometric fundamental group $\pi_1({\bf
G}_m\otimes {\bf F})$ and hence $H^i({\bf P}^1\otimes {\bf
F},j_\ast({\rm Sym}^k({\rm Kl}_n)))$ vanishes for $i=0,2$. So
$L({\bf P}^1, j_\ast({\rm Sym}^k({\rm Kl}_n)), T)={\rm det}(1-FT,
H^1({\bf P}^1\otimes {\bf F},j_\ast({\rm Sym}^k({\rm Kl}_n)))$ is
a polynomial. Since ${\rm Kl}_n$ is a lisse sheaf on ${\bf G}_m$
puncturely pure of weight $n-1$, all the reciprocal roots of the
polynomial ${\rm det}(1-FT, H^1({\bf P}^1\otimes {\bf
F},j_\ast({\rm Sym}^k({\rm Kl}_n)))$ have weight $k(n-1)+1$ by
[D2] 3.2.3. Similarly, one can show $L({\bf G}_m, {\rm Sym}^k({\rm
Kl}_n), T)$ is a polynomial. Let's prove it has coefficients in
${\bf Z}$. For any $a\in {\bf F}_p^*$, let $\sigma_a$ denote the
element of ${\rm Gal}({\bf Q}(\xi_p)/{\bf Q})$ such that
$\sigma(\xi_p) =\xi_p^a$, where $\xi_p\not=1$ is a $p$-th root of
unity in $\overline {\bf Q}$. Using the definition of Kloosterman
sums, one checks that
\[
\sigma_a({\rm Kl}_n({\bf F}_{q^k}, \lambda)) = {\rm Kl}_n({\bf
F}_{q^k}, a^n\lambda),
\]
and thus
\[
\sigma_a(L(\lambda, T)) = L(a^n\lambda, T).
\]
Recall from the Introduction that
\[
L(\lambda, T)^{(-1)^n} = (1-\pi_1(\lambda)T)\cdots
(1-\pi_n(\lambda)T).
\]
Write
\[
{\rm Sym}^k(L(\lambda, T)^{(-1)^n}) = \prod_{i_1+\cdots +i_n=k}
(1-\pi_1^{i_1}(\lambda)\cdots \pi_n^{i_n}(\lambda)T). \] Then, we
deduce
\[
\sigma_a({\rm Sym}^k(L(\lambda, T)^{(-1)^n})) = {\rm
Sym}^k(L(a^n\lambda, T)^{(-1)^n}). \] Now, by definition,
\[
L({\bf G}_m, {\rm Sym}^k({\rm Kl}_n), T) =\prod_{\lambda \in |{\bf
G}_m|} \frac{1}{ {\rm Sym}^k(L(\lambda, T^{{\rm
deg}(\lambda)})^{(-1)^n})}.\] Thus,
\[
\sigma_a(L({\bf G}_m, {\rm Sym}^k({\rm Kl}_n), T)) =\prod_{\lambda
\in |{\bf G}_m|} \frac{1}{ {\rm Sym}^k(L(a^n\lambda, T^{{\rm
deg}(\lambda)})^{(-1)^n})}.\] The right side is clearly the same
as $L({\bf G}_m, {\rm Sym}^k({\rm Kl}_n), T)$. The lemma is
proved.

\bigskip
\noindent {\bf Proposition 2.3.} Let $F_0$ be the geometric
Frobenius element at $0$ and $I_0$ the inertia subgroup at $0$.
The reciprocal roots of the polynomial ${\rm det}(1-F_0 T, ({\rm
Sym}^k({\rm Kl}_n))^{I_0})$ are of the form $q^i$, where $0\leq
i\leq k(n-1)/2$. In particular, ${\rm det}(1-F_0 T, ({\rm
Sym}^k({\rm Kl}_n))^{I_0})$ is a polynomial in integer
coefficients of weights at most $k(n-1)$. In the case $n=2$, we
have
$${\rm det}(1-F_0 T, ({\rm Sym}^k({\rm
Kl}_n))^{I_0})=1-T.$$

\bigskip
\noindent {\bf Proof.} By [K] 7.3.2 (3) and [D2] 1.8.1, the
eigenvalues of $F_0$ acting on ${\rm Sym}^k({\rm Kl}_n))^{I_0}$
are of the form $q^i$ with $0\leq i\leq k(n-1)/2$. This proves the
first part of the Proposition.

By [K] 7.4.3, the local monodromy at $0$ of ${\rm Kl}_n$ is
unipotent with a single Jordan block, and $F_0$ acts trivially on
$({\rm Kl}_n)^{I_0}$. If $n=2$, using this fact, one can show the
local monodromy at $0$ of ${\rm Sym}^k({\rm Kl}_2)$ has the same
property. Proposition 2.3 follows.

\bigskip
\noindent {\bf Remark}. We do not know a precise formula for ${\rm
det}(1-F_0 T, ({\rm Sym}^k({\rm Kl}_n))^{I_0})$ in general for
$n>2$.

\bigskip
\noindent {\bf Lemma 2.4.} Keep the notation in Lemma 1.7. Suppose
$pn$ is odd. Then $\theta$ is trivial.

\bigskip
\noindent {\bf Proof.} Take a positive integer $k$ such that
$(k,n)=1$ and such that $S_k(n,p)$ is non-empty. For instance, we
can take $k=n+mp$ for any positive integer $m$ prime to $n$.
(Recall that we always assume $(n,p)=1$). Then $a_k(n,p)\not= 0$.
By Lemma 2.2, $L({\bf G}_m,{\rm Sym}^k({\rm Kl}_n), T)$ is a
polynomial with coefficient in ${\bf Z}$. So each pure weight part
of $L({\bf G}_m,{\rm Sym}^k({\rm Kl}_n), T)$ also has coefficients
in ${\bf Z}$. We have
\begin{eqnarray*}
&&L({\bf G}_m,{\rm Sym}^k({\rm Kl}_n), T)\\
&=&L({\bf P}^1,j_\ast({\rm Sym}^k({\rm Kl}_n)), T) {\rm det}(1-F_0
T, ({\rm Sym}^k({\rm Kl}_n))^{I_0}){\rm det}(1-F_\infty T, ({\rm
Sym}^k({\rm Kl}_n))^{I_\infty}),
\end{eqnarray*}
where $F_\infty$ is the geometric Frobenius element at $\infty$.
Since $n$ is odd, by Lemmas 1.7 and 1.8, we have
$$({\rm Sym}^k({\rm Kl}_n))^{I_\infty} = ({\rm Sym}^k([n]_\ast{\cal L}_{\psi_n}
))^{I_\infty} \otimes {\cal L}_{\theta^k}\otimes \overline {\bf
Q}_l\left(\frac{k(1-n)}{2}\right).$$ Using the calculation in
Lemmas 1.8 and 2.1, one can verify $F_\infty$ acts trivially on
$({\rm Sym}^k ([n]_\ast{\cal L}_{\psi_n}))^{I_\infty}$. Let
$\lambda=\theta(F_\infty)$. Then $\lambda^n=1$ by Lemma 1.8 (1),
and
$${\rm det}(1-F_\infty T, ({\rm Sym}^k({\rm Kl}_n))^{I_\infty})
=(1-\lambda^k q^{\frac{k(n-1)}{2}}T)^{a_k(n,p)}.$$ So we have
$$L({\bf G}_m,{\rm Sym}^k({\rm Kl}_n), T)=$$
$${\rm det}(1-F_0T, ({\rm Sym}^k({\rm Kl}_n))^{I_0})
(1-\lambda^k q^{\frac{k(n-1)}{2}}T)^{a_k(n,p)}L({\bf
P}^1,j_\ast({\rm Sym}^k({\rm Kl}_n)), T). $$ By Lemma 2.2, $L({\bf
P}^1,j_\ast({\rm Sym}^k({\rm Kl}_n)), T)$ is pure of weight
$k(n-1)+1$. So the part of $L({\bf G}_m,{\rm Sym}^k({\rm Kl}_n),
T)$ with weight at most $k(n-1)$ is given by
$${\rm det}(1-F_0T, ({\rm Sym}^k({\rm Kl}_n))^{I_0})
(1-\lambda^kq^{\frac{k(n-1)}{2}}T)^{a_k(n,p)}.$$ It must have
coefficients in ${\bf Z}$. The first factor also has coefficients
in ${\bf Z}$. Working with the coefficients of $T$, we see that
$\lambda^k q^{\frac{k(n-1)}{2}}$ must be an integer. Since
$\lambda^n=1$, we must have $\lambda^k=\pm 1$. As $(k,n)=1$ and
$n$ is odd, we must have $\lambda=1$. So $\theta$ is trivial.

\bigskip
Now Theorem 1.1 in Section $1$ follows from Lemmas 1.7, 1.8, and
2.4.

\bigskip
In the following,  we calculate the bad factors at $\infty$ of the
$L$-function of the $k$-th symmetric product of ${\rm Kl}_n$.

\bigskip
\noindent{\bf Theorem 2.5.} Suppose $n|(q-1)$. Let $F_\infty$ be
the geometric Frobenius element at $\infty$.

(1) If $n$ is odd, then for all $k$, we have
\[{\rm det}(1-F_\infty T, ({\rm Sym}^k({\rm Kl}_n))^{I_\infty})=
(1-q^{\frac{k(n-1)}{2}}T)^{a_k(n,p)}.\]

(2) If $n$ is even and $k$ is odd, then we have
\[{\rm det}(1-F_\infty T, ({\rm Sym}^k({\rm Kl}_n))^{I_\infty})= 1.\]

(3) Suppose $n$ and $k$ are both even. We have
\begin{eqnarray*}
&&{\rm det}(1-F_\infty T, ({\rm Sym}^k({\rm Kl}_n))^{I_\infty})\\
&=& \left\{
\begin {array} {ll}
(1-q^{\frac{k(n-1)}{2}}T)^{b_k(n,p)}& {\rm if} \; 2n|(q-1), \\
(1+q^{\frac{k(n-1)}{2}}T)^{c_k(n,p)}(1-q^{\frac{k(n-1)}{2}}T)^
{b_k(n,p)-c_k(n,p)}& {\rm if} \; 2n\not|(q-1), \hbox { either }
4|n\;
{\rm  or}\; 4|k,\\
(1-q^{\frac{k(n-1)}{2}}T)^{c_k(n,p)}(1+q^{\frac{k(n-1)}{2}}T)^
{b_k(n,p)-c_k(n,p)}& {\rm if} \; 2n\not|(q-1), \;4\not|n \;{\rm
and}\; 4\not|k.
\end{array}\right.
\end{eqnarray*}
where $c_k(n,p)$ denotes the number of $\sigma$-orbits $j$ in
$S_k(n,p)$ such that $v_j\not=0$ and such that $j_1+2j_2+\cdots
+(n-1)j_{n-1}$ is odd.

\bigskip
\noindent {\bf Proof}. By Lemmas 1.7 and 1.8, we have
$$({\rm Sym}^k({\rm Kl}_n))^{I_\infty} = ({\rm Sym}^k([n]_\ast({\cal L}_{\psi_n}
\otimes{\cal L}_\chi)))^{I_\infty} \otimes {\cal
L}_{\theta^k}\otimes \overline {\bf
Q}_l\left(\frac{k(1-n)}{2}\right).$$ Suppose $n$ is odd. Then
$\chi$ and $\theta$ are trivial. One can verify $F_\infty$ acts
trivially on $({\rm Sym}^k ([n]_\ast{\cal
L}_{\psi_n}))^{I_\infty}$. (1) then follows from Lemma 2.1.

Suppose $n$ is even and $k$ is odd. Then $({\rm Sym}^k({\rm
Kl}_n))^{I_\infty}=0$ by Lemma 2.1. (2) follows.

Suppose $n$ and $k$ are even. If $2n|(q-1)$, then
$(\sqrt\zeta)^{q-1}=\zeta^{\frac{q-1}{2}}=1$ and hence $\sqrt
\zeta\in {\bf F}_q$. In this case, one verifies that $\theta$ is
trivial and $F_\infty$ acts trivially on $({\rm
Sym}^k([n]_\ast({\cal L}_{\psi_n} \otimes{\cal
L}_\chi)))^{I_\infty}$. The first case of (3) then follows from
Lemma 2.1. Suppose $2n\not|(q-1)$, then $\sqrt\zeta\not\in {\bf
F}_q$. We use the notation in the proof of Lemma 1.8. Note that
$g_{(0,1,-1)}$ is a lifting of the geometric Frobenius element in
${\rm Gal}({\bf F}/{\bf F}_p)$. Recall that $({\rm
Sym}^kV)^{I_\infty}$ is generated by the vectors
$$v_j = \sum_{i=0}^{n-1} (-1)^{j_{n-1}+\cdots +j_{n-i}} e^{\sigma^i(j)},$$
where $j$ runs over the $\sigma$-orbits of $S_k(n,p)$. One checks
that
\begin{eqnarray*}
&&g_{(0,1,-1)}(v_j)\\
&=&\sum_{i=0}^{n-1} (-1)^{j_{n-1}+\cdots +j_{n-i}} (-1)^{0\cdot
j_{n-i}+1\cdot j_{n-i+1}+\cdots +(i-1)j_{n-1}+ij_0+(i+1)j_1+\cdots
+(n-1)j_{n-i-1}} e^{\sigma^i(j)}.
\end{eqnarray*}
We have
\begin{eqnarray*}
&&0\cdot j_{n-i}+1\cdot j_{n-i+1}+\cdots+
(i-1)j_{n-1}+ij_0+(i+1)j_1+\cdots +(n-1)j_{n-i-1}\\
&=& i(j_0+j_1+\cdots +j_{n-i-1}+j_{n-i}+\cdots +j_{n-1})\\
&&+j_1+2j_2+\cdots
+(n-i-1)j_{n-i-1}\\
&&+(0-i)j_{n-i}+(1-i)j_{n-i+1}+\cdots
+(-1)j_{n-1}\\
&=& ik\\&&+\biggl(j_1+2j_2+\cdots
+(n-i-1)j_{n-i-1}\\&&+(n-i)j_{n-i}+(n-i+1)j_{n-i+1}+\cdots
+(n-1)j_{n-1}\biggr)
\\&&-n(j_{n-i}+j_{n-i+1}+\cdots +j_{n-1}).
\end{eqnarray*}
But $n$ and $k$ are even. So we  have
$$g_{(0,1,-1)}(v_j)=(-1)^{j_1+2j_2+\cdots +(n-1)j_{n-1}} v_j.$$
So $v_j$ is an eigenvector of $F_\infty$ on $({\rm
Sym}^kV)^{I_\infty}$ with eigenvalue $(-1)^{j_1+2j_2+\cdots
+(n-1)j_{n-1}}$. On the other hand, we have $\theta^k(F_\infty)=1$
if either $4|n$ or $4|k$, and $\theta^k(F_\infty)=-1$ if $4\not |
n$ and $4\not |k$. The last two cases of (3) follows.

\bigskip
\bigskip
\centerline {\bf 3. The degrees of the $L$-functions of symmetric
products}

\bigskip
\bigskip
Let $\zeta$ be a primitive $n$-th root of unity in ${\bf F}$ and
let $k$ be a positive integer. Denote by $d_k(n,p)$ the number of
the set $S_k(n,p)$, that is, the number of $n$-tuples
$(j_0,j_1,\ldots, j_{n-1})$ of non-negative integers satisfying
$j_0+j_1+\cdots + j_{n-1}=k$ and $j_0+j_1\zeta+\cdots + j_{n-1}
\zeta^{n-1}=0$ in ${\bf F}$. In this section, we prove the
following result, which is Theorem 0.1 in the Introduction:

\bigskip
\noindent {\bf Theorem 3.1.} Suppose $(n,p)=1$. The degree of
$L({\bf G}_m, {\rm Sym}^k({\rm Kl}_n), T)$ is
$$\frac{1}{n}\left({k+n-1\choose n-1} - d_k(n,p)\right).$$

\bigskip
\noindent {\bf Proof.} By Grothendieck's formula for
$L$-functions, $L({\bf G}_m, {\rm Sym}^k({\rm Kl}_n), T)$ is a
rational function, and its degree is the negative of the Euler
characteristics
$$\chi_c({\bf G}_m\otimes{\bf F},{\rm Sym}^k({\rm Kl}_n))= \sum_{i=0}^2 (-1)^i
{\rm dim}H^i_c({\bf G}_m\otimes {\bf F}, {\rm Sym}^k({\rm
Kl}_n))$$ To calculate the Euler characteristic, we may replace
the ground field ${\bf F}_q$ by its finite extensions. So we may
assume $n|(q-1)$. Then $\zeta$ lies in ${\bf F}_q$. By Lemmas 1.2
and 1.4, $[n]^\ast {\rm Kl}_n$ is isomorphic to ${\cal
L}_{\psi_n}\oplus {\cal L}_{\psi_{n \zeta}}\oplus\cdots \oplus
{\cal L}_{\psi_{n\zeta^{n-1}}}$ as a representation of the wild
inertia subgroup $P_\infty$ at $\infty$, where for any $a\in {\bf
F}_q$, $\psi_{a}$ is the additive character
$\psi_{a}(x)=\psi(ax)$. So we have
$$[n]^\ast({\rm Sym}^k({\rm Kl}_n))\cong \bigoplus_{j_0+j_1+\cdots+j_{n-1}=k,
\; j_0,j_1,\ldots,j_{n-1}\geq 0} {\cal
L}_{\psi_{n(j_0+j_1\zeta+\cdots+j_{n-1}\zeta^{n-1})}}$$ as
representations of $P_\infty$. But the Swan conductor of ${\cal
L}_{\psi_a}$ at $\infty$ is $1$ if $a\not =0$, and $0$ if $a=0$.
So the Swan conductor of $[n]^\ast ({\rm Sym}^k({\rm Kl}_n))$ at
$\infty$ is the number of those $n$-tuples $(j_0,j_1,\ldots,
j_{n-1})$ of non-negative integers satisfying $j_0+j_1+\cdots
+j_{n-1}=k$ and $j_0+j_1\zeta +\cdots +j_{n-1}\zeta^{n-1}\not=0$
in ${\bf F}_q$. This number is exactly ${k+n-1\choose n-1}
-d_k(n,p)$. By [K] 1.13.1, the Swan conductor of ${\rm Sym}^k({\rm
Kl}_n)$ at $\infty$ is $\frac{1}{n}$ of the Swan conductor of
$[n]^\ast ({\rm Sym}^k({\rm Kl}_n))$ at $\infty$. Since ${\rm
Kl}_n$ is tame at $0$, ${\rm Sym}^k({\rm Kl}_n)$ is also tame at
$0$, and hence its Swan conductor at $0$ vanishes. By the
Grothendieck-Ogg-Shafarevich formula, $-\chi({\bf G}_m,{\rm
Sym}^k({\rm Kl}_n))$ is equal to the sum of the Swan conductors of
${\rm Sym}^k ({\rm Kl}_n)$ at $0$ and at $\infty$. Theorem 3.1
follows.

\bigskip
In some special cases, an explicit formula for the degree of
$L({\bf G}_m, {\rm Sym}^k({\rm Kl}_n), T)$ can be obtained. The
following is one example.

\bigskip
\noindent {\bf Corollary 3.2.} Suppose $n$ is a prime number
different from $p$ such that $p$ is a primitive $(n-1)$-th root of
unity mod $n$. Let $\widetilde{(\frac{k}{n})}$ be the smallest
non-negative integer so that its image in ${\bf F}_p$ is
$\frac{k}{n}$. Then the degree of $L({\bf G}_m, {\rm Sym}^k({\rm
Kl}_n), T)$ is
$$\frac{1}{n}\left(\left(\begin{array}{c}
k+n-1 \\
n-1
\end{array}
\right)-\left(\begin{array}{c}
\frac{k-n \widetilde{({\frac{k}{n}})}}{p}+n-1 \\
n-1
\end{array}
\right)\right).$$

\bigskip
\noindent {\bf Proof.} Let $d=[{\bf F}_p(\zeta):{\bf F}_p]$. Then
we have $\zeta^{p^d-1}=1$. Since $\zeta$ has order $n$, we must
have $n|p^d-1$. On the other hand, since $p$ is a primitive
$(n-1)$-th root of unity mod $n$, $n-1$ is the smallest natural
number with the property $p^{n-1}=1$ in ${\bf Z}/n$, that is,
$n|p^{n-1}-1$. So we have $n-1\leq d$. Since $1+\zeta+\cdots +
\zeta^{n-1}=0$, we have $d=[{\bf F}_p(\zeta):{\bf F}_p]\leq n-1$.
So we must have $d=[{\bf F}_p(\zeta):{\bf F}_p]=n-1$ and
$1+X+\cdots + X^{n-1}$ is the minimal polynomial of $\zeta$ over
${\bf F}_p$. Therefore, if $j_0+j_1\zeta+ \cdots
+j_{n-1}\zeta^{n-1}=0$ in ${\bf F}$ for some integers
$j_0,j_1,\ldots, j_{n-1}$, then we must have
$$j_0\equiv j_1\equiv\cdots \equiv j_{n-1} \;({\rm mod} \; p).$$

Let us determine the number $d_k(n,p)$ of $n$-tuples
$(j_0,j_1,\ldots, j_{n-1})$ of non-negative integers with the
property $j_0+j_1+\cdots + j_{n-1}=k$ and $j_0+j_1\zeta +\cdots +
j_{n-1} \zeta^{n-1} =0$ in ${\bf F}$. By the above discussion, the
second equation implies that $j_0\equiv j_1\equiv\cdots \equiv
j_{n-1} \;({\rm mod} \; p)$. Substituting this into the first
equation, we get that
$$j_0\equiv j_1 \equiv\cdots \equiv j_{n-1} \equiv \widetilde{({\frac{k}{n}})}\;({\rm mod}
\; p).$$ Write $j_i=pj_i^\prime + \widetilde{({\frac{k}{n}})} $,
where $j_i^\prime$ are non-negative integers. Then the first
equation becomes $p(j_0^\prime+j_1^\prime+\cdots + j_{n-1}^\prime)
+ n\widetilde{({\frac{k}{n}})} =k,$ that is,
$$j_0^\prime+j_1^\prime+\cdots
+ j_{n-1}^\prime=\frac{k-n\widetilde{({\frac{k}{n}})}}{p}.$$ Thus
we have
$$d_k(n,p)=\left(\begin{array}{c}
\frac{k-n\widetilde{({\frac{k}{n}})}}{p}+n-1 \\ n-1
\end{array}
\right).$$ The corollary then follows from Theorem 3.1.

\bigskip
\bigskip
\centerline{\bf 4. An example}

\bigskip
\bigskip
In this section, we prove Theorem 0.2 in the Introduction.
Throughout this section, we assume that $n=2$, $q=p$ and $p$ is an
odd prime.

\bigskip
\noindent {\bf Lemma 4.1.} The degree of $L({\bf G}_m, {\rm
Sym}^k({\rm Kl}_2), T)$ is $\frac{k}{2}-[\frac{k}{2p}]$ if $k$ is
even, and $\frac{k+1}{2}-[\frac{k}{2p}+\frac{1}{2}]$ if $k$ is
odd.

\bigskip
\noindent {\bf Proof.} Note that $-1$ is a primitive square root.
Let's determine the number $d_k(2,p)$ of pairs $(j_0,j_1)$ of
non-negative integers satisfying $j_0+j_1=k$ and $j_0-j_1=0$ in
${\bf F}$, or equivalently,
\begin{eqnarray*}
j_0+j_1&=& k, \\
j_0-j_1&\equiv 0& (\hbox { mod } p).
\end{eqnarray*}
This is equivalent to the problem of determining the number of
those integers $0\leq j_0\leq k$ with the property $$2j_0\equiv k
\hbox { (mod } p).$$ Since $p$ is an odd prime, the inverse of $2$
in ${\bf F}_p$ is $\frac{p+1}{2}$, and hence the above equation is
equivalent to the equation $$j_0\equiv \frac{k(p+1)}{2} \hbox {
(mod } p).$$ So we are reduced to determining the number of those
integers $j$ with the property
$$0\leq jp+ \frac{k(p+1)}{2}\leq k,$$
that is,
$$-\frac{k}{2p}\leq j+\frac{k}{2}\leq \frac{k}{2p},$$
or equivalently,
$$-\frac{k}{2p}+\frac{1}{2}\leq j+\frac{k+1}{2}\leq
\frac{k}{2p}+\frac{1}{2}.$$ When $k$ is even, the number of those
integers $j$ satisfying $$-\frac{k}{2p}\leq j+\frac{k}{2}\leq
\frac{k}{2p}$$ is $2[\frac{k}{2p}]+1$. When $k$ is odd, the number
of those integers $j$ satisfying $$-\frac{k}{2p}+\frac{1}{2}\leq
j+\frac{k+1}{2}\leq \frac{k}{2p}+\frac{1}{2}$$ is
$2[\frac{k}{2p}+\frac {1}{2}]$. So
\[d_k(2,p)=\left\{\begin{array}{cc}
2[\frac{k}{2p}]+1& \hbox { if } k \hbox { is even},\\
2[\frac{k}{2p}+\frac{1}{2}]&\hbox { if } k \hbox { is odd}.
\end{array}
\right.\] The lemma then follows from Theorem 3.1.

\bigskip
\noindent {\bf Lemma 4.2.} Let $F_\infty$ be the geometric
Frobenius element at $\infty$. We have
\[{\rm det}(1-F_\infty T, ({\rm Sym}^k({\rm Kl}_2))^{I_\infty})=\left\{ \begin
{array} {ll}
1& {\rm if}\; 2\not | k, \\
(1-p^{\frac{k}{2}}T)^{m_k}& {\rm if} \; 2|k \;{\rm
and} \; p\equiv 1 \;({\rm mod}\; 4), \\
(1+p^{\frac{k}{2}}T)^{n_k}(1-p^{\frac{k}{2}}T)^ {m_k-n_k}& {\rm
if} \; 2|k \;{\rm and} \; p\equiv -1 \;({\rm mod}\; 4),
\end{array}
\right.\] where
\[m_k=\left\{ \begin {array} {ll}
1+[\frac{k}{2p}] & {\rm if}\; 4 | k, \\
{}[\frac{k}{2p}]& {\rm if} \; 4\not |k,
\end{array}
\right.\] and $n_k=[\frac{k}{4p}+\frac{1}{2}]$,

\bigskip
\noindent {\bf Proof.} The first case follows from Theorem 2.5
(2). Suppose $2|k$. We will treat the case where $p\equiv -1$ (mod
$4$) and leave to the reader to treat the other case. Note that
the condition $p\equiv -1$ (mod $4$) is equivalent to saying that
the square root $\sqrt {-1}$ of the primitive square root of unity
$-1$ does not lie in ${\bf F}_p$. We use the notation in the proof
of Lemma 1.8. The proof of Lemma 2.1 shows that a basis for $({\rm
Sym}^kV)^{I_\infty}$ is given by
$$\{e_0^ie_1^{k-i}+(-1)^{k-i}e_0^{k-i}e_1^i|i-(k-i)\equiv 0\;({\rm
mod}\; p),\; 0\leq i\leq \frac{k}{2}\}.$$ (When $i=\frac{k}{2}$,
the element
$e_0^ie_1^{k-i}+(-1)^{k-i}e_0^{k-i}e_1^i=(1+(-1)^{\frac{k}{2}})
e_0^{\frac{k}{2}}e_1^{\frac{k}{2}}$ is nonzero only when $4|k$. We
exclude this element from the basis if $4\not|k$.) A calculation
similar to that in the proof of Lemma 4.1 shows that this basis
has $1+[\frac{k}{2p}]$ elements if $4 | k$, and $[\frac{k}{2p}]$
elements if $4\not |k$. So we have
\[{\rm dim} ({\rm Sym}^k ({\rm
Kl}_2))^{I_\infty}=\left\{ \begin {array} {ll}
1+[\frac{k}{2p}] & {\rm if}\; 4 | k, \\
{}[\frac{k}{2p}]& {\rm if} \; 4\not |k.
\end{array}
\right.\]
Note that $g_{(0,1,-1)}$ is a lifting of the geometric
Frobenius element in ${\rm Gal}({\bf F}/{\bf F}_p)$. Let
$e_0^ie_1^{k-i}+(-1)^{k-i}e_0^{k-i}e_1^i$ be an element in the
above basis. Then
\begin{eqnarray*}
g_{(0,1,-1)}(e_0^ie_1^{k-i}+(-1)^{k-i}e_0^{k-i}e_1^i)& =&
(-1)^{k-i} e_0^ie_1^{k-i}+(-1)^{k-i}(-1)^ie_0^{k-i}e_1^i\\
&=&(-1)^i (e_0^ie_1^{k-i}+(-1)^{k-i}e_0^{k-i}e_1^i).
\end{eqnarray*}
(Recall that $k$ is even.) So $F_\infty$ acts semisimply on $({\rm
Sym}^kV)^{I_\infty}$ with eigenvalues $1$ and $-1$, and the
dimension of the eigenspace corresponding to the eigenvalue $1$
(resp. $-1$) is the number of those even (resp. odd) integers $i$
satisfying $0\leq i\leq \frac{k}{2}$ and $i-(k-i)\equiv 0$ (mod
$p$). (If $4\not |k$, we don't count the odd number
$i=\frac{k}{2}$). We have $i\equiv \frac{k}{2}$ (mod $p$). So
dimension of the eigenspace corresponding to the eigenvalue $1$
(resp. $-1$) is the number of integers $j$ such that $0\leq jp+
\frac{k}{2}\leq \frac{k}{2}$ and that $jp+\frac{k}{2}$ is even
(resp. odd). (Again if $4\not |k$, we don't count the odd number
$i=0\cdot p +\frac{k}{2}$). Note that if $4\not | k$, then
$jp+\frac{k}{2}$ is even if and only if $j$ is odd; if $4|k$, then
$jp+\frac{k}{2}$ is odd if and only if $j$ is odd. One can show
the number of odd integers $j$ satisfying $0\leq jp+
\frac{k}{2}\leq \frac{k}{2}$ is $[\frac{k}{4p}+\frac{1}{2}]$. So
if $4\not |k$ (resp. $4|k$), then the dimension of the eigenspace
corresponding to the eigenvalue $1$ (resp. $-1$) is
$[\frac{k}{4p}+\frac{1}{2}]$.

On the other hand, by Lemma 1.8, $({\rm Sym}^k({\rm
Kl}_2))^{I_\infty}$ is isomorphic to $({\rm Sym}^k
V)^{I_\infty}\otimes{\cal L}_{\theta^k}\otimes \overline {\bf
Q}_l\left(-\frac{k}{2}\right)$ as a representation of ${\rm
Gal}({\bf F}/{\bf F}_p)$. Since $k$ is even, we have
$\theta^k=(\theta^2)^{\frac{k}{2}}$. By the description of
$\theta^2$ in Lemma 1.8, if $4|k$, then $\theta^k(F_\infty)=1$; if
$4\not | k$, then $\theta^k(F_\infty)=-1$. Combining with the
above calculation, this proves Lemma 4.2.

\bigskip
Finally, we discuss the functional equation for $L({\bf G}_m, {\rm
Sym}^k({\rm Kl}_2), T)$. By Lemma 2.2, $L({\bf G}_m, {\rm
Sym}^k({\rm Kl}_2), T)$ and $L({\bf P}^1, j_\ast({\rm Sym}^k({\rm
Kl}_2)),T)$ are polynomials, where $j:{\bf G}_m\hookrightarrow
{\bf P}^1$ is the open immersion. By Proposition 2.3, we have
$${\rm det}(1-F_0T, ({\rm Sym}^k({\rm Kl}_2))^{I_0})=1-T.$$
Combining with Lemma 4.2, we see that $${\rm det}(1-F_0 T, ({\rm
Sym}^k({\rm Kl}_2))^{I_0}){\rm det}(1-F_\infty T, ({\rm
Sym}^k({\rm Kl}_2))^{I_\infty})$$ is exactly the polynomial $P_k$
defined in the Introduction. We have
\begin{eqnarray*}
&&L({\bf G}_m, {\rm Sym}^k({\rm Kl}_2), T)\\
&=&L({\bf P}^1, j_\ast({\rm Sym}^k({\rm Kl}_2)),T){\rm det}(1-F_0
T, ({\rm Sym}^k({\rm Kl}_2))^{I_0}){\rm det}(1-F_\infty T, ({\rm
Sym}^k({\rm Kl}_2))^{I_\infty}).
\end{eqnarray*}
So $L({\bf P}^1, j_\ast({\rm Sym}^k({\rm Kl}_2)),T)$ is the
polynomial $M_k$ defined in the Introduction.

By [K] 4.1.11, we have $({\rm Kl}_2)^\vee={\rm Kl}_2\otimes
\overline {\bf Q}_l(1)$. So $({\rm Sym}^k({\rm Kl}_2))^\vee={\rm
Sym}^k({\rm Kl}_2)\otimes \overline {\bf Q}_l(k)$ General theory
then shows that we have a functional equation
$$L({\bf P}^1, j_\ast({\rm Sym}^k({\rm Kl}_2)),T)=ct^\delta L({\bf
P}^1, j_\ast({\rm Sym}^k({\rm Kl}_2)),\frac{1}{p^{k+1}T}),$$ where
$$c=\prod_{i=0}^2{\rm det}(-F,H^i({\bf
P}^1\otimes {\bf F}, j_\ast({\rm Sym}^k({\rm Kl}_2))
)^{(-1)^{i+1}}$$ and $\delta=-\chi({\bf P}^1\otimes {\bf F},
j_\ast({\rm Sym}^k({\rm Kl}_2),T)$. This proves the functional
equation for $M_k$ in the Introduction.

\bigskip
\bigskip
\noindent {\bf References}

\bigskip
\bigskip
\noindent [D1] P. Deligne, {\it Applications de la Formule des
Traces aux Sommes Trigonom\'etriques}, in Cohomologie \'Etale (SGA
$4\frac{1}{2}$), 168-232, Lecture Notes in Math. 569,
Springer-Verlag 1977.

\bigskip
\noindent [D2] P. Deligne, {\it La Conjecture de Weil II}, Publ.
Math. IHES 52 (1980), 137-252.

\bigskip
\noindent [GK] E. Grosse-Kl\"onne, On families of pure slope
L-functions, {\it Documenta Math.}, 8(2003), 1-42.

\bigskip
\noindent [K] N. Katz, {\it Gauss Sums, Kloosterman Sums, and
Monodromy Groups}, Princeton University Press 1988.

\bigskip
\noindent [R] P. Robba, {\it Symmetric powers of $p$-adic Bessel
equation}, J. Reine Angew. Math. 366 (1986), 194-220.

\bigskip
\noindent [S] J.-P. Serre, {\it Linear Representations of Finite
Groups}, Sringer-Verlag 1977.

\bigskip
\noindent [SGA 4$\frac{1}{2}$] P. Deligne, et al, {\it Cohomologie
\'Etale} (SGA 4$\frac{1}{2}$), Lecture Notes in Math. 569,
Springer-Verlag 1977.

\bigskip
\noindent [Sp] S. Sperber, Congruence properties of
hyperkloosterman sums, Compositio Math., 40(1980), 3-33.

\bigskip
\noindent [ST] J.-P. Serre and J. Tate, {\it Good reduction of
abelian varieties}, Ann. of Math. 88 (1968), 492-517.

\bigskip
\noindent [W1] D. Wan, {\it Dimension variation of classical and
$p$-adic modular forms}, Invent. Math.,  133(1998), 449-463.

\bigskip
\noindent [W2] D. Wan, {\it Rank one case of Dwork's conjecture},
J. Amer. Math. Soc., 13(2000), 853-908.

\bigskip
\noindent [W3] D. Wan, {\it A quick introduction to Dwork's
conjecture}, Contemporary Math., Volume 245 (1999), 147-163.

\bigskip
\noindent [W4] D. Wan, {\it Geometric moment zeta functions}, in
Geometric Aspects of Dwork Theory, Walter de Gruyter,  2004,
1113-129.
\end{document}